\documentclass[MSNbibl,number,citesort,seceqn,dvips]{arxbj}

\aid{0}
\volume{19}
\issue{5A}
\pubyear{2013}
\firstpage{1750}
\lastpage{1775}
\doi{10.3150/12-BEJ428} 

\makeatletter

\newcommand{\filt}{({\mathcal{F}}_t)}
\newcommand{\filtwh}{(\wh{\mathcal{F}}_t)}
\newcommand{\filtn}{({\mathcal{F}}^n_t)}
\newcommand{\spa}{(\Omega, {\mathcal{F}}, ({\mathcal{F}}_t), P)}
\newcommand{\spawh}{(\wh\Omega, \wh{\mathcal{F}}, (\wh{\mathcal{F}}_t), \wh P)}
\newcommand{\R}{\mathbb{R}}
\newcommand{\cdw}{{\mathcal{C}}^2}
\newcommand{\dist}{\mbox{\rm dist}}
\newcommand{\Rp}{\mathbb{R}^+}
\newcommand{\Rd}{\mathbb{R}^d}
\newcommand{\D}{{\mathbb{D}} (\Rp,\Rd)}
\newcommand{\DD}{{\mathbb{D}} (\Rp,\R^{2d})}
\newcommand{\DDD}{{\mathbb{D}} (\Rp,\R^{3d})}
\newcommand{\Db}{{\mathbb{D}} (\Rp,\R^{d+1})}
\newcommand{\Daaa}{{\mathbb{D}} (\Rp,\R^{2d+1})}
\newcommand{\Dii}{\mathbb{D} (\Rp,\R^{2d})}
\newcommand{\Djj}{\mathbb{D} (\Rp,\R^{2})}

\newcommand{\Diii}{\mathbb{D} (\Rp,\R^{3d})}
\newcommand{\Cl}{\int_0^t}
\newcommand{\arrowss}{\mathop{\longrightarrow}\limits_{S}^{*}}
\newcommand{\lra}{\longrightarrow}
\newcommand{\ox}{\bar X^n}
\newcommand{\p}{{ P}}
\newcommand{\fn}{({\mathcal{F}}^n_t)}
\newcommand{\wh}{\widehat}
\newcommand{\pn}{P^n}
\newcommand{\rn}{^{\rho^n}}
\newcommand{\rnn}{\rho^n}
\newcommand{\tk}{t_{n,k}}
\newcommand{\tkk}{t_{n,k-1}}
\newcommand{\arrowd}{\mathop{\longrightarrow}\limits_{\mathcal{D}}}
\newcommand{\arrowds}{\mathop{\longrightarrow}\limits_{\mathcal{D}}^{*}}
\newcommand{\arrowdps}{\mathop{\longrightarrow}\limits_{\mathcal{P}}^{*}}
\newcommand{\arrowdq}{\mathop{\longrightarrow}_{{\mathcal{D}}({\mathcal{Q}}_{B})}}
\newcommand{\arrowp}{\mathop{\longrightarrow}\limits_{\mathcal{P}}}
\newcommand{\N}{\mathbb{N}}
\newcommand{\No}{\mathbb{N}\cup\{0\}}

\newremark{remark}{Remark}[section]
\newtheorem{theorem}[remark]{Theorem}
\newtheorem{proposition}[remark]{Proposition}
\newtheorem{lemma}[remark]{Lemma}
\newtheorem{corollary}[remark]{Corollary}
\newremark{example}[remark]{Example}
\makeatother

\begin{document}
\begin{frontmatter}

\title{Penalization methods for the Skorokhod problem and reflecting
SDEs with jumps}
\runtitle{Penalization methods}

\begin{aug}
\author{\fnms{Weronika} \snm{\L aukajtys}\thanksref{e1}\ead[label=e1,mark]{vera@mat.uni.torun.pl}} \and
\author{\fnms{Leszek} \snm{S\l omi\'nski}\corref{}\thanksref{e2}\ead[label=e2,mark]{leszeks@mat.uni.torun.pl}}
\runauthor{W. \L aukajtys and L. S\l omi\'nski} 
\address{Faculty of Mathematics and Computer Science, Nicolaus
Copernicus University, ul. Chopina 12/18 87--100 Toru\'n, Poland.
\printead{e1}; \printead*{e2}}
\end{aug}

\received{\smonth{5} \syear{2011}}
\revised{\smonth{10} \syear{2011}}

%
\begin{abstract}
We study the problem of approximation of solutions of the
Skorokhod problem and reflecting stochastic differential equations
(SDEs) with jumps by sequences of solutions of equations with
penalization terms. Applications to discrete approximation of weak
and strong solutions of reflecting SDEs are given. Our proofs are
based on new estimates for solutions of equations with
penalization terms and the theory of convergence in the Jakubowski
$S$-topology.
\end{abstract}

%

\end{frontmatter}
%

\section{Introduction}
Let $D$ be a convex open set in
$\R^d$. Consider a $d$-dimensional reflecting stochastic differential
equation (SDE),
\begin{equation}
\label{eq1.1} X_t=H_t+\int_0^tf(X_{s-})
 \,\mathrm{d}Z_s+K_t, \qquad t\in\Rp,
\end{equation}
where $Z$ is a $d$-dimensional semimartingale with $Z_0=0$, $H$ is
an adapted process with $H_0\in\bar D=D\cup\partial D$, and
$f\dvtx\R^d
\rightarrow\Rd\otimes\Rd$ is a continuous function such that
\begin{equation}
\label{eq1.2} \big\Vert f(x)\big\Vert \leq L\bigl(1+|x|\bigr),\qquad  x\in\Rd
\end{equation}
(for the precise definition, see Section \ref{sec3}).
Our main purpose is to study the problem of approximation of weak
solution of (\ref{eq1.1}) by solutions of nonreflecting SDEs of the form
\begin{equation}
\label{eq1.3} X^n_t=H^n_t+\int
_0^tf \bigl(X^n_{s-}
\bigr)  \,\mathrm{d}Z^n_s-n\int_0^t
\bigl(X^n_s-\Pi \bigl(X^n_s
\bigr) \bigr) \,\mathrm{d}s,  \qquad t\in\Rp, n\in\N,
\end{equation}
where $H^n$ and $Z^n$ are perturbations of $H$ and $Z$,
respectively, and $\Pi(x)$ denotes projection of $x$ on $\bar
D$. Because for large $n\in\N$, the drift term
$-n\int_0^t(X^n_s-\Pi(X^n_s))\, \mathrm{d}s$ forces $X^n$ to stay close
to~$\bar D$, it is called the penalization term, and the SDE
(\ref{eq1.3}) is called the SDE with penalization term.

The foregoing problem was intensively investigated in the case where
$f$ is a Lipschitz continuous function, $H=X_0$, and $Z$ is a
continuous semimartingale. In particular, Lions \textit{et al.} \cite{lms} and
Menaldi \cite{me} have proven that
$E\sup_{t\leq q}|X^n_t-X_t|^2\lra0$ for $q\in\Rp$, provided that
$Z$ is a $d$-dimensional standard Wiener process. In the case
where $Z$ has jumps, to the best of our knowledge, such a problem has
been considered previously only
by Menaldi and Robin \cite{mr} and \L aukajtys and S\l
omi\'nski \cite{laus}. Menaldi and Robin studied the case
where $Z$ is a diffusion with Poissonian jumps and $H=X_0$. However, they
imposed a very restrictive condition on the
Poissonian measure coefficient, and consequently, $K$ is a
process with continuous trajectories. In this case, earlier methods
of approximation remain in force. In earlier work, we
considered in detail the case where $H=X_0$ and $Z$ is a
general semimartingale.
Because the approximating sequence $\{X^n\}$ might not be relatively
compact in the Skorokhod topology $J_1$, we proved our
convergence results in the $S$-topology introduced by Jakubowski
\cite{ja}. It is worth pointing out that in both of the aforementioned
papers, the initial process $H$ is constant (i.e., $H=X_0$), and $f$ is
a Lipschitz continuous function.

The purpose of the present paper is to investigate the problem of
approximation of $X$ by $\{X^n\}$ in the case of arbitrary
initial process $H$ and arbitrary continuous coefficient $f$
satisfying the linear growth condition (\ref{eq1.2}). Our proofs
are based on new estimates for solutions of equations with
penalization terms.

The paper is organized as follows. In Section \ref{sec2} we consider a
deterministic problem of approximating a solution of the
Skorokhod problem $x_t=y_t+k_t$, $t\in\Rp$ on domain $D$
associated with a given function $y\in\D$ such that $y_0\in\bar
D$ (for precise definition, see Section \ref{sec2}). The penalization method
involves approximating $x$
by solutions
of equations of the form
\begin{equation}
\label{eq1.4} x^{n}_t=y^n_t-n\int
_0^t \bigl(x^{n}_s-\Pi
\bigl(x^{n}_s \bigr) \bigr) \,\mathrm{d}s, \qquad  t\in\Rp,
\end{equation}
where $y^n\lra y$ in the Skorokhod topology $J_1$. Lions and
Sznitman \cite{ls} and C\'epa \cite{ce} proved that $\{
x^n \}$ tends to $x$ if $y$ is continuous. We omit the latter
assumption and consider arbitrary function $y\in\D$. In this
general case, we prove that the variation of the penalization term
of the SDE (\ref{eq1.4}) is locally uniformly bounded and for
fixed $t\in\Rp$ $x^n_t\lra x_t$, provided that $\Delta y_t=0$,
which implies in particular that $x^n$ tends toward $x$ in the
$S$-topology. It is noteworthy that, similar to \cite{ce},
here we do not assume that the domain $D$ satisfies the so-called
condition ($\beta$) introduced by Tanaka \cite{ta}.

In Section \ref{sec3} we present new estimates on solutions $X^n$ of
equations with penalization terms associated with a given process
$Y^n$ such that $Y^n_0\in\bar D$, that is, solutions of SDEs,
\begin{equation}
\label{eq1.5} X^{n}_t=Y^n_t-n\int
_0^t \bigl(X^{n}_s-\Pi
\bigl(X^{n}_s \bigr) \bigr) \,\mathrm{d}s,\qquad   t\in\Rp.
\end{equation}
In particular,
we prove that if $Y^n$ is a process admitting the decomposition
$Y^n=H^n+M^n+V^n$, where $H^n$ is an $({\mathcal{F}}^n_t)$ adapted
process, $M^n$ is an $({\mathcal{F}}^n_t)$ adapted local martingale
with $M^n_0=0$ and $V^n$ is an $({\mathcal{F}}^n_t)$ adapted
processes of bounded variation with $V^n_0=0$, then, for every
$\delta, q>0$, $a\in D$ there exist constants $C_1,C_2>0$ such
that for every $\eta>0$,
\begin{eqnarray*}
P \Bigl(\sup_{t\leq q}\big|X^n_t-a\big|\geq\eta \Bigr)&
\leq& P \bigl(\omega'_{H^n}(\delta,q)\geq
d_a/2 \bigr)+P \Bigl(\sup_{t\leq
q}\big|H^n_t-a\big|
\geq C_1\eta \Bigr)
\\
&&{} +C_2\eta^{-2}E \bigl( \bigl[M^n
\bigr]_{q}+\big|V^n\big|_{q}^{2} \bigr)
\end{eqnarray*}
and
\begin{eqnarray*}
P \biggl(n\int_0^{q}\big|X^n_s-
\Pi \bigl(X^n_s \bigr)\big| \,\mathrm{d}s\geq
\eta^2 \biggr)& \leq&P \bigl(\omega'_{H^n}(\delta,q)
\geq d_a/2 \bigr)+7P \Bigl(\sup_{t\leq
q}\big|H^n_t-a\big|
\geq C_1\eta \Bigr)
\\
&&{} +C_2\eta^{-2}E \bigl( \bigl[M^n
\bigr]_{q}+\big|V^n\big|_{q}^{2} \bigr),
\end{eqnarray*}
where $\omega'$ denotes the usual modulus of continuity and
$d_a=\dist(a,\partial D)$.

In Section \ref{sec4}, we use estimates derived in Section \ref{sec3} to prove our main
results on the approximation of $X$ by $\{X^n\}$. We assume that
$\{Z^n\}$ is a sequence of semimartingales satisfying the
so-called condition (UT), and we prove that if $(H^n,Z^n)$
converges weakly to $(H,Z)$ in the $J_1$ topology,
then $\{X^n\}$ converges weakly in the $S$ topology to $X$. Moreover,
we prove convergence of finite-dimensional distributions of
$\{X^n\}$ to the corresponding finite-dimensional distributions of
$X$ outside the set of discontinuity points of $Z$ and $H$.
Consequently, using discrete approximations constructed in a manner
analogous to Euler's formula, we prove
the existence of a weak solution of the SDE (\ref{eq1.1}),
provided that $f$ is continuous and satisfies (\ref{eq1.2}). Moreover,
if the SDE (\ref{eq1.1}) has the weak uniqueness property,
then our approximations computed by simple recurrent
formulas allows us to obtain numerical solution of the SDE (\ref
{eq1.1}). In the case of reflected diffusion processes,
similar approximation schemes have been considered previously (see,
e.g., Liu \cite{li}, Pettersson \cite{p2}, S\l omi\'nski \cite{s4}).
In this section we also present some natural
conditions ensuring convergence of $\{X^n\}$ to $X$ in
probability provided that (\ref{eq1.1}) has the so-called
pathwise uniqueness property. Related results concerning diffusion
processes have been given by, for instance, Kaneko and Nakao
\cite{kn}, Gy\"ongy and Krylov \cite{gk}, Bahlali, Mezerdi and
Ouknine \cite{bmo}, Alibert and Bahlali \cite{ab} and S\l
omi\'nski \cite{s4}.

We note that we consider the space $\D$ equipped with two
different topologies, $J_1$~and~$S$. Definitions and required
results for the Skorokhod topology $J_1$ have been given by, for
example, Billingsley \cite{bil} and Jacod and Shiryayev
\cite{js}. For the convenience of the reader, we have collected basic
definitions and properties of the $S$-topology in the \hyperref[app]{Appendix}.
More details have been provided in Jakubowski \cite{ja}.

In this paper, we
use the following notation. Every process $X$ appearing in
the sequel is assumed to have trajectories in the space $\D$.
If $X=(X^1,\dots, X^d)$ is a semimartingale, then $[X]_t$ represents
$\sum^d_{i=1}[X^i]_t$ and
$[X^i]$ represents the quadratic variation process of $X^i$, $i=1,\dots
, d$.
Similarly, $\langle X\rangle_t
=\sum^d_{i=1}\langle X^i\rangle_t$, and
$\langle X^i\rangle$ represents the predictable compensator
of $[X^i]$, $i=1,\dots, d$. If $K=(K^1,\dots, K^d)$ is the process
with locally finite
variation, then $|K|_t=\sum^d_{i=1}|K^i|_t$,
where $|K^i|_t$ is a total
variation of $K^i$ on $[0, t]$. In general,\vadjust{\goodbreak} we let ${
\arrowd}$ and ${ \arrowp}$ denote
convergence in law and in probability, respectively. To avoid
ambiguity,
we write ${  X^n\arrowd X}$ (${
X^n\arrowp X}$) in $\D$ if $\{X^n\}$ converges weakly (in probability)
to $X$
in the space $\D$ equipped with $J_1$.
Following \cite{ja}, we write ${  X^n\arrowds X}$
(${
X^n\arrowdps X}$) in $\D$ when we consider the $S$ topology.
For $x \in
\D$, $\delta>0$, $q\in\Rp$, we let $\omega'_x(\delta, q)$
and $\omega''_x(\delta, q)$ denote classical moduli of continuity of $x$
on $[0,q]$, that is, $\omega'_x(\delta, q)=\inf\{ \max_{i \leq r}
\omega_x([t_{i-1},t_i)); 0=t_0 < \cdots< t_r =q, \inf_{i <
r}(t_i - t_{i-1}) \geq\delta\}$, and ${\omega}''_{x}(\delta, q)
= \sup\{ \min(|x_u - x_s|,|x_t-x_u|), 0 \leq s < u < t \leq q, t
- s < \delta\},$
where $\omega_x(I)=\sup_{s,t \in I}|x_s-x_t|$. We
also use the modulus $\bar{\omega}''_{(x,y)}(\delta, q)$
introduced in Jakubowski \cite{ja}. We recall that for $x,y \in
\D$, $\delta
> 0$, $q \in\Rp$
$\bar{\omega}''_{(x,y)}(\delta, q) = \sup\{ \min(|x_u -
x_s|,|y_t-y_u|), 0 \leq s < u < t \leq q, t - s < \delta\}$.

\section{A deterministic case} \label{sec2}
Let $D$ be a
nonempty convex (possibly unbounded) open set in $\mathbb{R}^d$,
and let
${\mathcal{N}}_x$ denote the set of
inward normal unit vectors at $x\in\partial D$ ($\mathbf{n}\in{\mathcal{N}}_{x} $ if and only if ${\langle}y-x,\mathbf{n}{\rangle}\geq0
$ for every $y\in\bar D$, where
${\langle} \cdot , \cdot {\rangle}$ denotes the usual inner
product in $\mathbb{R}^d$). The following remark also can be found in
Menaldi \cite{me} or Storm \cite{sto}.

\begin{remark}\label{rem2.1}
(i) If $\dist(x,\bar D)>0$, then there exists a
unique $\Pi(x)\in\partial D$ such that $|\Pi(x)-x|=\dist(x,\bar
D)$. Moreover,
$(\Pi(x)-x)/|\Pi(x)-x|\in{\mathcal{N}}_{\Pi(x)}$.

(ii) For every $a\in D$,
\[
\big|x-\Pi(x)\big|\leq d_a^{-1}{ \bigl\langle}x-a,x-\Pi(x){ \bigr
\rangle}, \qquad x\in \Rd,
\]
where $d_a=\dist(a,{\partial
D})$.
\end{remark}

Let $y\in\D$ be a function with initial value in $\bar D$. We
recall that a pair of functions
$(x, k)\in\Dii$ is called a solution of the
Skorokhod problem associated with $y$ if
\begin{itemize}
\item[$\bullet$]   $x_t=y_t+k_t,\  t\in\Rp$,

\item[$\bullet$]  $x$ is $\bar D$-valued,

\item[$\bullet$]  $k$ is a function with locally bounded variation such that
$k_0=0$ and
\[
  k_t=\int_0^t\mathbf{n}_s \,\mathrm{d}|k|_s,\qquad |k|_t=\int
_0^t\mathbf{1}_{\{x_s\in\partial D\}} \,\mathrm{d}|k|_s,
 \qquad t\in\mathbb R^+,
\]
where $\mathbf{n}_s\in{\mathcal{N}}_{x_s}$ if
$x_s\in\partial D$.
\end{itemize}
The problem of existence of solutions of the Skorokhod problem and
its approximation by solutions of equations with penalization
terms has been discussed by many authors. Tanaka \cite{ta} proved
existence and uniqueness of solutions in the case of continuous
$y$ and domains also satisfying the following condition:
\begin{enumerate}[$(\beta)$]
\item[$(\beta)$] there exist constants $\epsilon>0$ and
$\delta>0$ such that for every $x\in\partial D$, we can find
$x_0\in D$ such that $B(x_0,\epsilon)\subset D$ and
$|x-x_0|\leq\delta$.
\end{enumerate}\eject
Tanaka also observed that $(\beta)$ holds true in dimensions 1
and 2 or if $D$ is a bounded set. On the other hand, in
dimension $>$2, one can construct examples of
nonbounded convex domains not satisfying ($\beta$). For
instance, the cone with the basis $\{(x,y,0)\in\R^3; y\geq
x^2\}$ and peak at $(0,0,1)$, that is, the set
\begin{equation}
\label{eq2.1}
C= \bigl\{(\lambda x,\lambda y,1-\lambda)\in\R^3;
x^2\leq y, 0\leq \lambda\leq1 \bigr\},
\end{equation}
does not satisfy ($\beta$). C\'epa \cite{ce} omitted the
assumption $(\beta)$ and proved the existence and uniqueness of
the solution to the Skorokhod problem in the case of continuous
function $y$. In addition, C\'epa proved convergence $\sup_{t\leq
q}|x^n_t-x_t|\lra0$, $q\in\Rp$ of solutions of equations
(\ref{eq1.4}) for every sequence $\{y^n\}$ such that $\sup_{t\leq
q}|y^n_t-y_t|\lra0$, $q\in\Rp$.

The case of functions
with jumps was considered for the first time by Anulova and
Liptser~\cite{al}, who proved the existence and uniqueness of
solutions under condition ($\beta$). Their result was generalized
to the case of arbitrary convex $D$ by \L aukajtys \cite{lau}.
In an earlier work \cite{laus}, we considered the problem of
approximating noncontinuous $x$ by solutions of
equations with penalization terms only in this very special case.
We now consider the problem
of approximating noncontinuous $x$ by solutions of equations
with penalization terms in the general case of arbitrary sequences
$\{y^n\}$ such that $y^n\lra y$ in $\D$. Our main tools are the
following estimates on the solution of (\ref{eq1.4}):

\begin{lemma}\label{lm2}
Let $y^n\in\D$, $y^n_0\in\bar{D}$ and let $x^n$ be
a solution of the equation (\ref{eq1.4}). Then for any $q,\delta
>0$ and $a \in D$ such that
\begin{equation}
\label{eq2.3}
\omega'_{y^n}(\delta, q) <
d_a/{2},
\end{equation}
we have
\begin{longlist}
\item[(i)]${
\sup_{t \leq q} |x^n_t-a| \leq2 \sqrt{7}([q/\delta]+1)
\sup_{t\leq q}|y^n_t-a|,}$
\item[(ii)]${
|k^n|_q \leq 55([q/ \delta]+1)^3d_a^{-1}\sup_{t\leq
q}|y^n_t-a|^2,}$
\end{longlist}
where $k^n_t=-n\int_0^t(x^{n}_s-\Pi(x^{n}_s)) \,\mathrm{d}s$, $ t\in\Rp$ and
$[q/\delta]$ denotes the largest integer less or equal to~$q/\delta$.
\end{lemma}

\begin{pf}
We follow the proof of Theorem 3.2 in \cite{ce}. Let $0 \leq t
\leq q $. Because $k^n$ is a continuous function such that
$k^n_0=0$,
\begin{eqnarray*}
\big|x^n_t-a\big|^2 & = &\big|y^n_t-a\big|^2
+ \bigl\langle k^n_t,k^n_t \bigr
\rangle+2 \Cl \bigl\langle y^n_t-a, \mathrm{d}k^n_u
\bigr\rangle
\\[-2pt]
& = & \big|y^n_t-a\big|^2+2\Cl
\bigl\langle x^n_u-a,\mathrm{d}k^n_u
\bigr\rangle+2\Cl \bigl\langle y^n_t-y^n_u,\mathrm{d}k^n_u
\bigr\rangle.
\end{eqnarray*}
Therefore, for any $0 \leq s \leq t \leq q$,
\begin{eqnarray*}
\big|x^n_t-a\big|^2-\big|x^n_s-a\big|^2
& = & \big|y^n_t-a\big|^2 - \big|y^n_s-a\big|^2
+2\int_s^t \bigl\langle x^n_u-a,\mathrm{d}k^n_u
\bigr\rangle
\\[-2pt]
& &{} - 2\int_s^t \bigl\langle
y^n_u-y^n_s,\mathrm{d}k^n_u
\bigr\rangle +2 \bigl\langle k^n_t,y^n_t-y^n_s
\bigr\rangle.\vadjust{\goodbreak}
\end{eqnarray*}
By Remark \ref{rem2.1}(ii),
\begin{eqnarray*}
2 \int_s^t \bigl\langle x^n_u
- a , \mathrm{d} k^n_u \bigr\rangle& = & -2n\int
_s^t \bigl\langle x^n_u
- a ,x^n_u- \Pi \bigl(x^n_u
\bigr) \bigr\rangle\, \mathrm{d}u
\\
& \leq& -2 d_a n \int_s^t
\big|x^n_u - \Pi \bigl(x^n_u
\bigr)\big|\,\mathrm{d}u = -2d_a\big|k^n\big|_s^t,
\end{eqnarray*}
and, consequently,
\begin{eqnarray*}
\big|x^n_t-a\big|^2-\big|x^n_s-a\big|^2
& \leq& \big|y^n_t-a\big|^2 - \big|y^n_s-a\big|^2-
2 d_a\big|k^n\big|_s^t - 2\int
_s^t \bigl\langle y^n_u-y^n_s,\mathrm{d}k^n_u
\bigr\rangle
\\
& &{} -2 \bigl\langle y^n_t-a,y^n_t-y^n_s
\bigr\rangle- 2 \bigl\langle a-x^n_t,y^n_t-y^n_s
\bigr\rangle
\\
& \leq& 5\sup_{t\leq q}\big|y^n_t-a\big|^2 + 4
\sup_{t\leq q}\big|y^n_t-a\big|\cdot \sup_{t\leq q}\big|x^n_t-a\big|
-2d_a\big|k^n\big|^t_s
\\
& &{} -2 \int_s^t \bigl\langle
y^n_u-y^n_s,\mathrm{d}k^n_u
\bigr\rangle.
\end{eqnarray*}
By (\ref{eq2.3}), there exists a subdivision $(s_k)$ of $[0,q]$
such that $0=s_0 < s_1 < \cdots<s_r = q$, $\delta\leq s_k -
s_{k-1}$, $k=1, \ldots,r-1$, where $r=[q/\delta]+1$ and
$\omega_{y^n}([s_{k-1},s_k))<d_a/2$. Thus, in particular,
\[
\int_{s_{k-1}}^{s_k} \bigl\langle y^n_u-y^n_{s_{k-1}},\mathrm{d}k_u
\bigr\rangle \leq\bigg| \int_{(s_{k-1},s_{k})} \bigl\langle
y^n_u-y^n_{s_{k-1}},\mathrm{d}k^n_u
\bigr\rangle\bigg|\leq\frac
{d_a}{2}\big|k^n\big|_{s_{k-1}}^{s_k}.
\]
Therefore,
\[
2 \biggl(-\int_{s_{k-1}}^{s_k} \bigl\langle
y^n_u-y^n_{s_{k-1}},\mathrm{d}k^n_u
\bigr\rangle - d_a\big|k^n\big|_{s_{k-1}}^{s_k}
\biggr) \leq 2 \biggl(\frac{d_a}{2}\big|k^n\big|_{s_{k-1}}^{s_k}
- d_a\big|k^n\big|_{s_{k-1}}^{s_k} \biggr) =
-d_a\big|k^n\big|_{s_{k-1}}^{s_k},
\]
which implies that
\begin{equation}
\label{eq2.4}
\big|x^n_{s_k}-a\big|^2-\big|x^n_{s_{k-1}}-a\big|^2
\leq 5\sup_{t\leq q}\big|y^n_t-a\big|^2 + 4
\sup_{t\leq q}\big|y^n_t-a\big|\cdot \sup_{t\leq q}\big|x^n_t-a\big|
-d_a\big|k^n\big|_{s_{k-1}}^{s_k}.
\end{equation}
From (\ref{eq2.4}), it follows immediately that
\begin{eqnarray*}
\big|x^n_{s_k}-a\big|^2-\big|x^n_{s_{k-1}}-a\big|^2
& \leq&5\sup_{t\leq
q}\big|y^n_t-a\big|^2 +4
\sup_{t\leq q}\big|y^n_t-a\big|\cdot\sup_{t \leq
q}\big|x^n_t-a\big|.
\end{eqnarray*}
Set $k_0=\max\{ k, s_k \leq t \}$. Then
\begin{eqnarray*}
\big|x^n_t-a\big|^2 &=& \sum
_{k=1}^{k_0} \bigl(\big|x^n_{s_k}-a\big|^2-\big|x^n_{s_{k-1}}-a\big|^2
\bigr)
\\
&& {}  +\big|x^n_t-a\big|^2-\big|x^n_{s_{k_0}}-a\big|^2+\big|x^n_0-a\big|^2
\\
&   \leq& r \Bigl(5\sup_{t\leq q}\big|y^n_t-a\big|^2+4
\sup_{t\leq
q}\big|y^n_t-a\big|\cdot \sup_{t \leq q}\big|x^n_t-a\big|
\Bigr)+\sup_{t\leq
q}\big|y^n_t-a\big|^2,
\end{eqnarray*}
which implies that
\[
\sup_{t \leq q}\big|x^n_t-a\big|^2 \leq
14r^2\sup_{t\leq q}\big|y^n_t-a\big|^2+
\sup_{t \leq q}\big|x^n_t-a\big|^2/2.
\]
Thus, $\sup_{t \leq q }|x^n_t-a|^2 \leq28r^2\sup_{t\leq q}|y^n_t-a|^2$,
and the proof of (i) is complete.

(ii) Using (\ref{eq2.4}) and (i) gives
\begin{eqnarray*}
d_a\big|k^n\big|_{s_{k-1}}^{s_k} &\leq&5
\sup_{t\leq
q}\big|y^n_t-a\big|^2+4
\sup_{t\leq q}\big|y^n_t-a\big|\cdot\sup_{t \leq q
}\big|x^n_t-a\big|
\\
&& {}     +\big|x^n_{s_{k-1}}-a\big|^2-\big|x^n_{s_k}-a\big|^2
\\
&   \leq& 13\sup_{t\leq q}\big|y^n_t-a\big|^2+
\frac{3}{2}\sup_{t \leq q}\big|x^n_t-a\big|^2
\leq55r^2\sup_{t\leq q}\big|y^n_t-a\big|^2
\end{eqnarray*}
for $k=1,\ldots,r$. Thus, $ |k^n|_q \leq
\sum_{k=1}^{r}|k^n|_{s_{k-1}}^{s_k}
\leq 55r^3d_a^{-1}\sup_{t\leq q}|y^n_t-a|^2$,
which completes the proof.
\end{pf}

\begin{theorem}
\label{tw2.15} Assume that $\{y^n\}\subset\D$, $y^n_0 \in
\bar{D}$, and let $x^n$ denote the solution of the equation
(\ref{eq1.4}), $n\in\N$. If $y^n\lra y$ in $\D$, then
\begin{longlist}
\item[(i)] $  \sup_{n \in\N}\sup_{t \leq
q}|x^n_t|< +\infty, q\in\Rp$ and $\sup_{n \in\N}|k^n|_q < +\infty$, $q\in\Rp$,
\item[(ii)]
$ x^n_t \lra x_t$,  provided that  $|\Delta y_t|=0$,
\item[(iii)]$ {x^n {\arrowss} x}$,
\end{longlist}
where $x$ denotes the solution of the Skorokhod problem
associated with $y$.
\end{theorem}

\begin{pf}
(i) Because $\{y^n\}$ is relatively compact in $\D$,
$\sup_n\sup_{t\leq q} |y^n_t|<\infty$ for any \mbox{$q\in\Rp$},
and for any $a\in D$,
$q>0$, there exists $\delta>0$ such that
$\sup_n\omega'_{y^n}(\delta,q)<d_a/2$. Therefore, the first
conclusion follows from Lemma \ref{lm2}.

(ii) Let $\{ \delta^i \}$ be a sequence of constants such that
$\delta^i \downarrow0$ and $|\Delta y_t| \neq\delta^i , t \in
\Rp$. Set \mbox{$t^i_{n,0}= 0$}, $ t^i_{n,k+1}=\min( t^i_{n,k} +
\delta_k^i, \inf\{ t
> t^i_{n,k},|\Delta y^n_t|> \delta^i \})$,
$t^i_{0}=0, $ $
t^i_{k+1}=\min( t^i_{k} + \delta_k^i, \inf\{ t >
t^i_{k}, |\Delta y_t|> \delta^i \})$, where $\{\{ \delta^i_k
\}\}$ is an array of constants satisfying $\delta^i /2 \leq
\delta_k^i \leq\delta^i$ and $|\Delta y_{t^i_k + \delta
_k^i}|=0$. Now, for every $i\in\N$, set
$ y^{n,(i)}_t = y^n_{t^i_{n,k}}$, $t \in
[t^i_{n,k},t^i_{n,k+1})$, $ y^{(i)}_t=y_{t^i_{k}}$, $t \in
[t^i_{k},t^i_{k+1})$, $n, k \in\No$. Observe that for every $q
\in\Rp$,
\[
\lim_{i\to\infty}\sup_{t\leq q}\big|y^{(i)}_t -
y_t\big|=0\quad  \mbox{and}\quad   \lim_{i\to\infty}
\limsup_{n \rightarrow\infty}\sup_{t\leq
q}\big|y^{n,(i)}_t -
y^n_t\big|=0 .
\]
Let $x^{n,(i)}$ be a solution of
an equation with a penalization term of the form
\[
x^{n,(i)}_t=y^{n,(i)}_t-n\int
_0^t \bigl(x^{n,(i)}_s-\Pi
\bigl(x^{n,(i)}_s \bigr) \bigr)  \,\mathrm{d}s=y^{n,(i)}_t
+ k^{n,(i)}_t , \qquad t \in\Rp.
\]
Fix $t \in\Rp$ and consider the
decomposition $ x^n_t - x_t = x^n_t - x^{n,(i)}_t+ x^{n,(i)}_t -
x^{(i)}_t+x^{(i)}_t-x_t $, where $x^{(i)}$ denotes a solution of
the Skorokhod problem associated with $y^{(i)}$. Due to
\cite{laus}, Lemma 2.2(i),
\[
\big|x^n_t - x^{n,(i)}_t\big|^2
\leq\big|y^n_t - y^{n,(i)}_t\big|^2
+ 4\sup_{s \leq
t}\big|y^n_s-y^{n,(i)}_s\big|
\bigl(\big|k^n\big|_t+\big|k^{n,(i)}\big|_t \bigr),
\]
where variations $|k^n|_t$, $|k^{n,(i)}|_t$ are bounded uniformly
by Lemma \ref{lm2}(ii). Therefore,
\[
\lim_{i \rightarrow\infty} \limsup_{n \rightarrow\infty} \big|x^n_t -
x^{n,(i)}_t\big|=0.
\]
Moreover, if $|\Delta y_t| = 0$, then ${ \lim_{i
\rightarrow\infty}|\Delta y^{(i)}_t|=0}$. Because, by \cite{laus}, Lemma
3.3,\linebreak[4] $\limsup_{n \rightarrow\infty}|x^{n,(i)}_t-
x^{(i)}_t| \leq|\Delta y^{(i)}_t|$, $i\in\N$, it follows that $
\lim_{i \rightarrow\infty}\limsup_{n \rightarrow\infty
}|x^{n,(i)}_t- x^{(i)}_t| \leq\lim_{i \rightarrow\infty}|\Delta
y^{(i)}_t|=0$. On the other hand, by \cite{ta}, Lemma 2.2,
$\sup_{t \leq q}|x^{(i)}_t - x_t|\lra0$, $q \in\Rp,$ and (ii)
follows.

(iii) The sequence $\{ y^n \}$ is relatively compact in $J_1$, and
consequently it is relatively $S$ compact. Because by part (i),
$\sup_{n \in\N}|k^n|_q < +\infty$, $q \in\Rp$, the sequence
$\{k^n\}$ is also relatively $S$ compact, and thus $\{x^n\}$ is
relatively $S$ compact as well. In view of Corollary \ref{ap1}
(\hyperref[app]{Appendix}), this proves~(iii).
\end{pf}

Recall that if $y^n \lra y$ in $\D$, then for every $t
\in\Rp$ there exists a sequence $t_n \lra t$ such that
\begin{equation}
\label{zbdet} y^n_{t_n} \lra y_t,\qquad
y^n_{t_n-} \lra y_{t-} \quad  \mbox{and}\quad
\Delta y^n_{t_n} \lra\Delta y_t.
\end{equation}
Moreover, for arbitrary sequences $\{t_n'\}$, $\{t_n''\}$ such
that $t_n' < t_n \leq t''_n$, $n\in\N$, and $\lim_{n \rightarrow
\infty} t'_n=\lim_{n \rightarrow\infty} t_n''=t$, we have
\begin{equation}
\label{zbprim} y^n_{t_n'} \lra y_{t-}
\quad  \mbox{and}\quad   y^n_{t_n''} \lra y_t
\end{equation}
(see, e.g., \cite{js}, Chapter VI, Proposition 2.1).

\begin{corollary}\label{tw3.11}
Under the assumptions of Theorem
\ref{tw2.15},
\begin{longlist}
\item[(i)]
for every $t \in\Rp$, if $t_n \lra t$ is a sequence satisfying
(\ref{zbdet}), then
\[
x^n_{t_n} \lra x_{t-}+\Delta y_t
\]
and for arbitrary sequences $\{t_n'\}$, $\{t_n''\}$ such that
$t_n' < t_n < t''_n$, $n\in\N$ and $\lim_{n \rightarrow\infty}
t'_n=\lim_{n \rightarrow\infty} t_n''=t$, we have
\[
x^n_{t_n'} \lra x_{t-}  \quad  \mbox{and}\quad
x^n_{t_n''} \lra x_t ,
\]
\item[(ii)]
$(\Pi(x^n),y^n) \lra(x,y) \mbox{ in } \DD,$
\item[(iii)] moreover, if $y$ is continuous, then
\[
\sup_{t\leq q}\big|x^n_t-x_t\big|\lra0,\qquad  q
\in\Rp.\vadjust{\goodbreak}
\]
\end{longlist}
\end{corollary}

\begin{pf}
(i) If $y$, $y^n$ are step functions, then the result follows from
\cite{laus}, Lemma 3.3. In the general case, it is sufficient to
use (\ref{zbdet}), (\ref{zbprim}) and repeat the approximation
procedure from the proof of Theorem \ref{tw2.15}.

(ii) By Theorem \ref{tw2.15}(ii), $\Pi(x^n_t) \lra x_t$, provided
that $\Delta y_t=0$. Therefore, it suffices to prove that $\{(\Pi
(x^n ) , y^n)\}$ is relatively compact in $\DD$. Because $\sup_{n
\in\N}\break\sup_{t \leq q}|\Pi(x^n_t)| < + \infty$, $\sup_{n \in
\N}\sup_{t \leq q}|y^n_t| < + \infty$, $q\in\Rp$, it is sufficient
to show that
\begingroup
\abovedisplayskip=7pt
\belowdisplayskip=7pt
\begin{equation}
\label{limsupRzut} \lim_{\delta\downarrow0} \limsup_{n \rightarrow
\infty}
\omega_{(\Pi(x^n),y^n)} \bigl([0,\delta] \bigr)=0
\end{equation}
and
\begin{equation}
\label{limomega} \lim_{\delta\downarrow0} \limsup_{n \rightarrow
\infty}
\omega''_{(\Pi(x^n ), y^n )}(\delta, q)=0, \qquad  q \in
\Rp.
\end{equation}
To prove (\ref{limsupRzut}), first observe that
$x^n_0 \lra x_0$ and $y^n_0\lra0$, which implies that
(\ref{limsupRzut}) is equivalent to the following:
\begin{itemize}
\item[$\bullet$]
for every sequence $\{s^n_0\}$ such that $0 \leq s^n_0 \lra0$,
\begin{equation}
\label{war3.1} \Pi \bigl(x^n_{s_0^n} \bigr) \lra
x_0, \qquad y^n_{s_0^n}\lra y_0.
\end{equation}
\end{itemize}
Next, note that (\ref{war3.1}) is implied by (i) (it is sufficient
to put $t=0$ and observe that in this case, $t_n=0$).

Similarly, (\ref{limomega}) is equivalent to the condition
\begin{itemize}
\item[$\bullet$] for every $t\leq q$ and every sequence
$\{s^n_i\}$, $i=1,2,3$ such that $s^n_1 \leq s^n_2 \leq
s^n_3$,  $n\in\N$, if $\lim_{n \rightarrow\infty}s^n_i=t$,
$\Pi(x^n_{s^n_i}) \rightarrow a_i$, $y^n_{s^n_i} \rightarrow
b_i$, $i=1,2,3$, then
\begin{equation}
\label{war3.2} a_1=a_2\quad  \mbox{and}\quad
b_1=b_2 \quad  \mbox{or}\quad   a_2=a_3
\quad  \mbox{and}\quad  b_2=b_3.
\end{equation}
\end{itemize}
Because $\Pi$ is continuous and $x_t=\Pi(x_{t-}+\Delta y_t)$, it
follows from (i) that for arbitrary sequences~$\{t_n'\}$,
$\{t_n''\}$ such that $t_n' < t_n \leq t''_n$ and $\lim_{n
\rightarrow\infty} t'_n=\lim_{n \rightarrow\infty} t''_n=t$, we
have
\begin{equation}
\label{zbrzut} \Pi \bigl(x^n_{t_n'} \bigr) \lra
x_{t-} \quad  \mbox{and}\quad  \Pi \bigl(x^n_{t_n''}
\bigr) \lra x_{t}.
\end{equation}
Combining (\ref{zbprim}) with (\ref{zbrzut}), we see that there are
only four possibilities:
\begin{eqnarray*}
 a_1&=&a_2=a_3=x_{t-}
\quad  \mbox{and}\quad b_1=b_2=b_3=y_{t-},
\\[-1pt]
a_1&=&a_2=x_{t-},\qquad   a_3=x_t
\quad  \mbox{and}\quad b_1=b_2=y_{t-},\qquad
b_3=y_{t},
\\[-1pt]
a_1&=&x_{t-},\qquad   a_2=a_3=x_{t}
\quad  \mbox{and}\quad b_1=y_{t-},\qquad   b_2=b_3=y_{t},
\\[-1pt]
a_1&=&a_2=a_3=x_{t} \quad  \mbox{and}\quad b_1=b_2=b_3=y_{t}.
\end{eqnarray*}
Thus, in each case (\ref{war3.2}) is satisfied.

(iii) In this case, $x$ is continuous as well. Moreover,
(\ref{zbdet}) is satisfied for every $t \in\Rp$ and every
sequence $t_n \lra t$. Consequently, by part (i), for every $t
\in\Rp$ and every sequence $t_n \lra t$,
\[
x^n_{t_n} \lra x_{t-}+\Delta
y_t=x_t,
\]
which is equivalent to (iii).
\end{pf}\endgroup\eject

\begin{corollary}
\label{wniosek3.14}
Let $\{ (y^n,z^n) \}$ be relatively compact in
$\DD$, $y^n_0 \in\bar{D}$, $n\in\N$. If $x^n$ denotes the
solution of (\ref{eq1.4}), $n\in\N$, then, for every $q \in\Rp$,
\[
\lim_{\delta\downarrow0} \limsup_{n \rightarrow\infty
}\bar{\omega}''_{(x^n,z^n)}(
\delta, q)=0.
\]
\end{corollary}

\begin{pf}
Without loss of generality, we may and will assume that $(y^n,z^n)
\lra(y,z)$ in $\DD$. Then, for every $t \in\Rp$, there is a
sequence $t_n \lra t$ such that
\begin{equation}
\label{zb3.15} y^n_{t'_n} \lra y_{t-},\qquad
 z^n_{t'_n} \lra z_{t-} \quad  \mbox{and}\quad
y^n_{t''_n} \lra y_{t}, \qquad  z^n_{t''_n}
\lra z_{t}
\end{equation}
for arbitrary sequences $t'_n < t_n \leq t''_n$ such that $\lim_{n
\rightarrow\infty} t'_n=\lim_{n \rightarrow\infty} t''_n=t$.
Because\linebreak[4] $ \sup_{n \in\N}\sup_{t \leq q}|x^n_t| < + \infty$
and $\sup_{n \in\N}\sup_{t \leq q}|y^n_t| < +\infty$, $q\in\Rp$,
the proof is completed by showing that
\begin{itemize}
\item[$\bullet$] for every $t \leq q$ and every sequence $\{s^n_i\}$ such that
$\lim_{n \rightarrow\infty}s^n_i=t$, $i=1,2,3$ and $s^n_1\leq\allowbreak
s^n_2 \leq s^n_3$, $n\in\N$, if there exists limits
$x^n_{s^n_1} \rightarrow a_1$, $x^n_{s^n_2} \rightarrow a_2$,
and
$ z^n_{s^n_2} \rightarrow b_1$, $ z^n_{s^n_3}
\rightarrow b_2$, then
\begin{equation}
\label{war3.3}a_1=a_2    \quad  \mbox{or}\quad
b_1=b_2 .
\end{equation}
\end{itemize}
From (\ref{zbprim}), (\ref{zb3.15}), and Theorem \ref{tw3.11}(i),
we conclude that there are only the following cases:
\begin{eqnarray*}
 a_1&=&a_2=x_{t-}
  \quad  \mbox{and}\quad b_1=b_2=z_{t-}  \quad  \mbox{or}\quad  b_1=z_{t-} \quad  \mbox{and}\quad b_2=z_t,
\\
a_1&=&x_{t-}, \qquad a_2=x_{t-}+\Delta
y_t   \quad  \mbox{and}\quad b_1=b_2=z_{t},
\\
a_1&=&x_{t-},\qquad  a_2=x_t \quad  \mbox{and}\quad b_1=b_2=z_t,
\\
a_1&=&a_2=x_{t-}+\Delta y_t  \quad  \mbox{and}\quad b_1=b_2=z_{t},
\\
a_1&=&x_{t-}+\Delta y_t, \qquad a_2=x_t
  \quad  \mbox{and}\quad b_1=b_2=z_{t},
\\
a_1&=&a_2=x_t  \quad  \mbox{and}\quad
b_1=b_2=z_t;
\end{eqnarray*}
that is, (\ref{war3.3}) is satisfied.
\end{pf}

\section{Applications to stochastic processes}\label{sec3}

Let $Y^n$ be an $\filtn$ adapted process with $Y^n_0 \in\bar{D}$
and let $X^n$ be
a solution of the equation (\ref{eq1.5}).
We will consider processes $Y^n$ admitting the decomposition
\begin{equation}
\label{eq3.1} Y^n_t=H^n_t+M^n_t+V^n_t,
 \qquad t\in\Rp,
\end{equation}
where $H^n$ is an $({\mathcal{F}}^n_t)$ adapted process,
$M^n$ is an $({\mathcal{F}}^n_t)$ adapted local martingale with
$M^n_0=0$, and $V^n$ is an $({\mathcal{F}}^n_t)$ adapted processes of
bounded variation with $V^n_0=0$.

\begin{remark}[(\cite{laus})]
\label{uwaga}
Let $\wh Y^n$ be another $\filtn$ adapted processes
of the form $\wh Y^n=H^n+\allowbreak  \wh M^n+\wh
V^n$, where $\wh M^n$ is a local martingale, $ \wh V^n$ is a
process with locally bounded variation, and \mbox{$\wh M^n_0=\wh
V^n_0=0$}. Assume that $\wh Y^n_0\in\bar D$, and let $\wh X^n$ be a
solution of the equation
\begin{equation}
\label{eq3.2} \wh X^n_t=\wh Y^n_t-n
\int_0^t \bigl(\wh X^n_s-
\Pi \bigl(\wh X^n_s \bigr) \bigr) \,\mathrm{d}s, \qquad  t\in\Rp.
\end{equation}
For every $p \in\mathbb{N}$, there exists $C(p)$ such that
\begin{equation}
\label{eq3.4}
E\sup_{t \leq\tau}\big|X^n_t-\wh
X^n_t\big|^{2p} \leq C(p) E \bigl(
\bigl[M^n-\wh M^n \bigr]_{\tau}^p +
\big|V^n-\wh V^n\big|_{\tau}^{2p} \bigr)
\end{equation}
and
\begin{equation}
\label{eq3.13} E\sup_{t < \tau}\big|X^n_t-\wh
X^n_t\big|^{2p} \leq C(p) E \bigl(
\bigl[M^n-\wh M^n \bigr]_{\tau-}^p +
\big|V^n-\wh V^n\big|_{\tau-}^{2p} + \bigl
\langle M^n-\wh M^n \bigr\rangle^p_{\tau-}
\bigr)
\end{equation}
for every stopping time $\tau$.
\end{remark}

\begin{theorem}\label{tw2}
Let $Y^n$ be a process of the form (\ref{eq3.1}), and let $X^n$
denote the solution of (\ref{eq1.5}). For
every $\delta, q$, $a\in D$, there exist constants $C_1,C_2>0$
such that for every $\eta>0$,\vspace*{-1pt}
\begin{eqnarray}\label{eq3.5}
P \Bigl(\sup_{t\leq q}\big|X^n_t-a\big|\geq\eta \Bigr)
& \leq&  P \bigl(\omega'_{H^n}(\delta,q)\geq
d_a/2 \bigr)+P \Bigl(\sup_{t\leq q}\big|H^n_t-a\big|
\geq C_1\eta \Bigr)
\nonumber\\[-8pt]\\[-8pt]
&& {}   +C_2\eta^{-2} E \bigl(
\bigl[M^n \bigr]_{q}+\big|V^n\big|_{q}^{2}
\bigr)
\nonumber
\end{eqnarray}
and\vspace*{-1pt}
\begin{eqnarray}\label{eq3.6}
\hspace*{-12pt}P \biggl(n\int_0^q\big|{X}_{s}^n-
\Pi \bigl({X}_{s}^n \bigr)\big|\,\mathrm{d}s\geq \eta^2
\biggr) &\leq& P \bigl(\omega'_{H^n}(\delta,q)
\geq d_a/2 \bigr) +7P \Bigl(\sup_{t\leq q}\big|H^n_t-a\big|
\geq C_1\eta \Bigr)\hspace*{12pt}
\nonumber\\[-8pt]\\[-8pt]
&&  {}  +C_2\eta^{-2}E \bigl(
\bigl[M^n \bigr]_{q}+\big|V^n\big|_{q}^{2}
\bigr).
\nonumber
\end{eqnarray}
\end{theorem}

\begin{pf}
Let $C'=2\sqrt{7}([T/\delta]+1)$ be a constant from
Lemma \ref{lm2}(i), and let $\wh X^n$ be a solution of equation
with penalization term (\ref{eq3.2}) associated with $\wh
Y^n=H^n$. Then\vspace*{-1pt}
\begin{eqnarray*}
\p \Bigl(\sup_{t \leq q}\big|\wh X^n_t-a\big|\geq\eta/2
\Bigr)& = & \p \Bigl(\sup_{t \leq
q}\big|\wh X^n_t-a\big|\geq
\eta/2 , \omega_{H^n}'(\delta,q) < d_a/2
\Bigr)
\\
& & {}+ \p \bigl(\omega_{H^n}'(\delta,q)\geq
d_a/2 \bigr)
\\
&\leq& \p \Bigl(\sup_{t \leq q}\big|H^n_t-a\big|
\geq{C'}^{-1} \eta/2 \Bigr) +\p \bigl(
\omega_{H^n}'( \delta,q)\geq d_a/2 \bigr).
\end{eqnarray*}
By the foregoing and by (\ref{eq3.4}),\vspace*{-1pt}
\begin{eqnarray*}
\p \Bigl(\sup_{t \leq q}\big|X^n_t-a\big|\geq\eta \Bigr) &
\leq& \p \Bigl( \sup_{t \leq
q}\big|\wh X^n_t-a\big| \geq
\eta/2 \Bigr)+\p \Bigl( \sup_{t \leq q}\big|X^n_t - \wh
X^n_t\big| \geq\eta/2 \Bigr)
\\
& \leq& \p \Bigl(\sup_{t \leq q}\big|H^n_t-a\big|\geq\eta
\bigl(2C' \bigr)^{-1} \Bigr) +\p \bigl(
\omega_{H^n}'( \delta,q) \geq d_a/2 \bigr)
\\
& &{} + 4C(1) \eta^{-2} E \bigl( \bigl[M^n
\bigr]_q + \big|V^n\big|^2_q \bigr),
\end{eqnarray*}
which completes the proof of (\ref{eq3.5}) with $C_1=2C'$ and
$C_2=4C(1)$.

Now, for simplicity of notation, set
$K_t^n=-n\int_0^t({X}_{s}^n-\Pi({X}_{s}^n))\,\mathrm{d}s$. Our proof of
(\ref{eq3.6}) starts with the observation that the estimates
similar to (\ref{eq3.5}) is true for $\sup_{t\leq q}|K^n_t|$ as well.
Indeed, there exist constants $C_1,C_2>0$ such that for every
$\eta>0$,
\begin{eqnarray}\label{eq3.7}
\p \Bigl(\sup_{t \leq q}\big|K^n_t\big|\geq \eta
\Bigr) & \leq& \p \Bigl(\sup_{t
\leq
q}\big|X^n_t-a\big|\geq
\eta/3 \Bigr) + \p \Bigl(\sup_{t \leq q}\big|H^n_t-a\big|\geq
\eta/3 \Bigr)\nonumber
\\
& &
{}+ \p \Bigl(\sup_{t \leq q}\big|M^n_t+V_t^n\big|
\geq\eta/3 \Bigr)
\nonumber\\[-8pt]\\[-8pt]
& \leq&  2\p \Bigl(\sup_{t \leq q}\big|H^n_t-a\big|>
C_1\eta \Bigr) +\p \bigl(\omega_{H^n}'(
\delta,q) \geq d_a/2 \bigr)\nonumber
\\
& & {}+ C_2 \eta^{-2}E \bigl( \bigl[M^n
\bigr]_q + \big|V^n\big|^2_q \bigr)
\nonumber
.
\end{eqnarray}
On the other hand, by Remark \ref{rem2.1}(ii),
\begin{eqnarray*}
\big|K^n\big|_q & \leq& -\frac{1}{d_a} \int
_0^q \bigl\langle X^n_s-a
, \mathrm{d}K^n_s \bigr\rangle
\\
& = &-\frac{1}{d_a} \biggl(\int_0^q \bigl
\langle H^n_s-a,\mathrm{d}K^n_s \bigr
\rangle+ \int_0^q \bigl\langle
Z^n_s , \mathrm{d}K^n_s \bigr\rangle+
\int_0^q \bigl\langle K^n_s
, \mathrm{d}K^n_s \bigr\rangle \biggr),
\end{eqnarray*}
where $Z^n=M^n+V^n$. Because $K^n$ has continuous trajectories, it
follows by the integration by parts formula that $\int_0^q \langle
Z^n_s , \mathrm{d}K^n_s \rangle = \langle Z^n_q , K^n_q \rangle-
\int_0^q \langle K^n_{s} , \mathrm{d}Z^n_s \rangle$, and $\int_0^q \langle
K^n_s , \mathrm{d}K^n_s \rangle = \frac{1}{2}|K_q|^2\geq0$. Therefore,
\begin{equation}
\label{eq3.8} \big|K^n\big|_q \leq-\frac{1}{d_a} \biggl(\int
_0^q \bigl\langle H^n_s-a,\mathrm{d}K^n_s
\bigr\rangle+ \bigl\langle Z^n_q , K^n_q
\bigr\rangle- \int_0^q \bigl\langle
K^n_{s} , \mathrm{d}Z^n_s \bigr\rangle
\biggr).
\end{equation}
Fix $\omega \in\{ \omega_{H^n }' (\delta,q)  <   d_a/2
 \}$. There exists a subdivision $(s_k)$ of $[0,q]$ such that
$\delta\leq s_k-s_{k-1}, k=1,2,\ldots,r-1$, where $r=[q/\delta]+1 $
and $\omega_{H^n(\omega)}([s_{k-1},s_k))<d_a/2$. Set
$ H_t^{n*}(\omega) = H^n_{s_{k-1}}(\omega), $ for $   t
\in[s_{k-1} ,s_{k}) $. Then $\sup_{t<q} |H_t^{n}(\omega
)-H_t^{n*}(\omega)|<d_a/2$ and, consequently,
\begin{eqnarray*}
\bigg|\int_0^q \bigl\langle H^n_s(
\omega)-a , \mathrm{d}K^n_s(\omega) \bigr\rangle \bigg|&\leq&\bigg|\int
_0^q \bigl\langle H^n_s(
\omega)-H^{n*}(\omega) , \mathrm{d}K^n_s(\omega) \bigr
\rangle\bigg|
\\
&&{} +\bigg|\int_0^q \bigl\langle
H^{n*}_s(\omega)-a , \mathrm{d}K^n_s(
\omega) \bigr\rangle\bigg|
\\
& \leq& \frac{d_a}{2}\big|K^n\big|_q(\omega) + 2r
\sup_{t \leq
q}\big|H^n_t(\omega)-a\big|\big| K^n_t(
\omega)\big|.
\end{eqnarray*}
Combining this inequality with (\ref{eq3.8}), we see that on the
set $\{ \omega_{H^n }' (\delta,q)  <   d_a/2  \}$, we have
\[
\big|K^n\big|_q \leq2d_a^{-1} \biggl( 2r
\sup_{t \leq q}\big|H^n_t-a\big|\big| K^n_t\big|
+ \sup_{t \leq
q}\big|Z^n_t\big|\big| K^n_t\big|
+ \bigg|\int_0^q \bigl\langle K^n_{s}
, \mathrm{d}Z^n_s \bigr\rangle\bigg| \biggr).
\]
Thus, there is a constant $C>0$ such that
\begin{eqnarray}\label{eq3.9}
P \bigl(\big|K^n\big|_q \geq\eta^2 \bigr) &\leq& P
\bigl(\big|K^n\big|_q \geq \eta^2,
\omega_{H^n}'(\delta,q) < d_a/2 \bigr) + P
\bigl(\omega_{H^n}'(\delta ,q)\geq d_a/2
\bigr)
\nonumber
\\
&\leq& P \Bigl( \sup_{t \leq
q}\big|H^n_t-a\big|\big| K^n_t\big|
\geq C \eta^2 \Bigr) + P \Bigl(\sup_{t \leq q}\big|Z^n_t\big|\big| K^n_t\big|
\geq C\eta^2 \Bigr)
\\
& &{} + P \biggl(\bigg|\int_0^q \bigl\langle
K^n_{s}, \mathrm{d}Z^n_s \bigr\rangle\bigg|
\geq C\eta^2 \biggr) + P \bigl(\omega_{H^n}'(
\delta,q) \geq d_a/2 \bigr).
\nonumber
\end{eqnarray}
Clearly,
\[
P \Bigl( \sup_{t \leq q}\big|H^n_t-a\big|\big| K^n_t\big|
\geq C \eta^2 \Bigr)\leq P \Bigl( \sup_{t \leq
q}\big|K^n_t\big|
\geq \eta \Bigr)+P \Bigl( \sup_{t \leq q}\big|H^n_t-a\big|
\geq C \eta \Bigr)
\]
and
\begin{eqnarray*}
P \Bigl(\sup_{t \leq q}\big|Z^n_t\big|\big| K^n_t\big|
\geq C\eta^2 \Bigr)&\leq& P \Bigl(\sup_{t
\leq q}\big|K^n_t\big|
\geq\eta \Bigr) +P \Bigl( \sup_{t \leq q}\big|Z^n_t\big|\geq
C \eta \Bigr)
\\
&\leq&P \Bigl(\sup_{t \leq q}\big|K^n_t\big| \geq\eta
\Bigr)+C_2 \eta^{-2} E \bigl( \bigl[M^n
\bigr]_q + \big|V^n\big|^2_q \bigr).
\end{eqnarray*}
Moreover, if we set $ \tau^n = \inf\{ t, |K^n_t| \geq\eta \} \wedge
q$, then, obviously, $\sup_{t \leq\tau^n }|K^n_t| \leq\eta$ and
\begin{eqnarray*}
\p \biggl(\bigg|\int_0^q \bigl\langle
K^n_{s} , \mathrm{d}Z^n_s \bigr\rangle\bigg|
\geq C \eta^2 \biggr) &\leq& \p \bigl( q > \tau^n \bigr) +
\p \biggl(\bigg|\int_0^{ \tau^n } \bigl\langle
K^n_{s} , \mathrm{d}Z^n_s \bigr\rangle\bigg|
\geq C \eta^2 \biggr)
\\
& \leq& \p \Bigl(\sup_{t \leq q }\big|K^n_t\big|\geq\eta
\Bigr) + C_2\eta^{-4} E \biggl(\int_0^{ \tau^n }
\bigl\langle K^n_{s} , \mathrm{d}Z^n_s
\bigr\rangle \biggr)^2
\\
& \leq& \p \Bigl(\sup_{t \leq q }\big|K^n_t\big|\geq\eta
\Bigr) + C_2 \eta^{-2} E \bigl( \bigl[M^n
\bigr]_q + \big|V^n\big|^2_q \bigr).
\end{eqnarray*}
Combining the last three inequalities with (\ref{eq3.7}) and
(\ref{eq3.9}) yields (\ref{eq3.6}).
\end{pf}

Now, for $n\in\N$, let $Z^n$ be a semimartingale adapted to some
filtration $({\mathcal{F}}^n_t)$. We assume that $\{Z^n\}$
satisfies the following condition (UT) introduced by Stricker
\cite{st}:
\begin{enumerate}[(UT)]
\item[(UT)] For every $q\in\Rp$, the
family of random variables
\[
\biggl\{\int_{[0,q]} U^n_s
 \,\mathrm{d}Z^n_s; n\in\N , U^n\in\mathbf{U}^n_q \biggr\}
\]
is bounded in probability.
Here $ \mathbf{U}^n_q$ is the class of discrete predictable
processes of the form $U^n_s=U^n_0+\sum_{i=0}^kU^n_i\mathbf{
1}_{\{t_i<s\leq t_{i+1}\}}$, where $0=t_0<t_1<\cdots<t_k=q$ and
$U^n_i$ is ${\mathcal{F}}^n_{t_i}$ measurable, $|U^n_i|\leq1$ for
$i\in\{0,\ldots,k\} , n,k\in\N.$
\end{enumerate}
The condition (UT) proved to be very useful in the theory of
limit theorems for stochastic integrals and for solutions of SDEs
(see, e.g., \cite{jmp,kp,ms,s1,s2}).

\begin{corollary}
\label{thm2.7} Let $\{Y^n\}$ be a sequence of $\filtn$ adapted
processes, and let $\{X^n\}$ be a sequence of solutions of
equations with penalization terms (\ref{eq1.5}). Assume that every
$Y^n$ is of the form $Y^n=H^n+Z^n$ with $H^n_0\in\bar D$,
$Z^n_0=0$, where $\{H^n\}$ is tight in the $\D$ sequence of $\filtn$
adapted processes and $\{Z^n\}$ is a sequence of $\filtn$-adapted
semimartingales satisfying \textup{(UT)}. Then, for every $q\in\Rp$, the
sequence $\{n\int_0^q|X^n_s-\Pi(X^n_s)| \,\mathrm{d}s\}$
is bounded in probability.
\end{corollary}

\begin{pf}
Define $\tau^n_k=\inf\{t;|H^n_t|\vee|Z^n_t|> k\}$, $k,n\in\N$.
Because (UT) implies that $\{\sup_{t\leq q}|Z^n_t|\}$ is bounded
in probability and $\{\sup_{t\leq q}|H^n_t|\}$ is bounded in
probability by the tightness of $ \{H^n\}$ in $\D$, we have
\begin{equation}
\label{eq3.10} \lim_{k\rightarrow\infty}\limsup_{n\rightarrow\infty}{ P} \bigl(
\tau^n_k\leq q \bigr)=0,\qquad  q\in\Rp.
\end{equation}
Furthermore, by simple calculations for every $\delta>0$ and $q\in\Rp$
$ \omega'_{H^{n,\tau^n_k-}}(\delta,q) \leq
\omega'_{H^{n}(\delta,q)}$\vspace*{-3pt} and $ \sup_{t\leq
q}|H^{n,\tau^n_k-}_t|\leq\sup_{t\leq q}|H^{n}_t|$, and thus the
sequence $\{H^{n,\tau^n_k-}\}$ also is tight in $\D$. On the other
hand, from the definition of (UT), we see that
$\{Z^{n,\tau^n_k}\}$ satisfies (UT) as well. Moreover, because
$ Z^{n,\tau^n_k-}_{\cdot} = Z^{n,\tau^n_k}_{\cdot}- \Delta
Z^n_{\tau^n_k}\mbox{\bf1}_{\{\cdot\geq\tau^n_k\}}$ and $|\Delta
Z^n_{\cdot}| \leq2|Z^n_{\cdot}|$, it follows from the definition
of (UT) that $\{Z^{n,\tau^n_k-}\}$ satisfies (UT) as well. Therefore,
in view of
(\ref{eq3.10}), without
loss of generality, we can and will assume that
$H^n=H^{n,\tau^n_k-}$ and $Z^n=Z^{n,\tau^n_k-}$ for some
$k\in\Rp$. Then $Z^n$ admits the decomposition
$Z^n=M^n+V^n$ with $|\Delta M^n|, |\Delta V^n|\leq4k$, where
$\{[M^n]_q\}$, $\{|V^n|_q\}$ are bounded in probability for each
$q\in\Rp$ (see,
e.g., \cite{ms}). Set
$\gamma^n_b=\inf\{t ;|H^n|_t\vee [M^n]_t\vee\allowbreak  |V^n|_t> b\}$ for
$n\in\N$, $b\in\Rp$. Then
\begin{equation}
\label{eq3.11} \lim_{b\rightarrow\infty}\limsup_{n\rightarrow\infty}{ P} \bigl(
\gamma^n_b\leq q \bigr)=0,\qquad  q\in\Rp,
\end{equation}
so as before, we can assume that $H^n=H^{n,\gamma^n_b}$,
$M^n=M^{n,\gamma^n_b} ,V^n=V^{n,\gamma^n_b} $ for some
$b\in\Rp$, and thus that $[M^n]_{\infty}\leq b+16k^2
 , |V^n|^2_{\infty}\leq b^2+16k^2$.

To complete the proof, it suffices to use Theorem
\ref{tw2}. Fix $\varepsilon>0$, $q\in\Rp$, and $ a\in D$, and let
$\delta>0$, $n_0\in\N$ be such that $P(\omega'_{H^n}(\delta,q)\geq
d_a/2)\leq\varepsilon/2$ for every $n\geq n_0$. If we put
$\eta=\max((k+|a|+1)/C_1,\sqrt{2C_2(b^2+b+32k^2)/\varepsilon})$ in
(\ref{eq3.6}) then
\[
P \biggl(n\int_0^q\big|X^n_s-
\Pi \bigl(X^n_s \bigr)\big|\, \mathrm{d}s\geq\eta^2
\biggr) \leq\varepsilon
\]
for
any $n\geq n_0$, which completes the proof.
\end{pf}

\begin{corollary}
\label{cor2.8} For $n\in\N$, let $Y^n$ and $\wh Y^n$ be
processes adapted to filtrations
$({\mathcal{F}}_t^n)$ and $(\wh{\mathcal{F}}_t^n)$, respectively, and let
$X^n$ be a solution of
(\ref{eq1.5}) and $\wh X^n$ be
a solution of (\ref{eq3.2}).
If\vspace*{1pt} $\{Y^n=H^n+Z^n\}$, $\{\wh Y^n=\wh H^n+\wh Z^n\}$ with
$H^n_0,\wh H^n_0\in\bar D$ and $Z^n_0=\wh Z^n_0=0$, and
$\{H^n\}$, $\{\wh H^n\}$ are tight in $\D$, $\{Z^n\}$, $\{\wh
Z^n\}$ satisfy \textup{(UT)} and
\[
\sup_{t\leq q}\big|Y^n_t-\wh Y^n_t\big|
\arrowp0,\qquad  q\in\Rp
\]
then
\[
\sup_{t\leq q}\big|X^n_t-\wh X^n_t\big|
\arrowp0,\qquad  q\in\Rp.
\]
\end{corollary}

\begin{pf}
The proof follows immediately from \cite{laus}, Lemma 2.2(i), and Corollary
\ref{thm2.7}.
\end{pf}

\begin{corollary}\label{wnzb} Let $\{Y^n\}$
be a given sequence of processes, $Y^n_0\in\bar D$, $n\in\N$, and
let $\{X^n\}$ be a sequence of solutions of equations with
penalization terms (\ref{eq1.5}):
\begin{longlist}
\item[(i)] For any sequence of processes $\{Z^n\}$, if
\[
\bigl\{ \bigl(Y^n,Z^n \bigr) \bigr\}\qquad \mbox{is tight
in }  {\mathbb{D}} \bigl(\Rp,\R^{2d}\bigr)
\]
then for
every $\varepsilon>0$, $q \in\Rp$
\begin{equation}
\label{eq3.12}
\lim_{\delta\downarrow0}\limsup_{n \rightarrow
\infty} \pn \bigl(\bar{
\omega}''_{(X^n,Z^n)}(\delta, q)> \varepsilon
\bigr)=0.
\end{equation}
\item[(ii)] For any sequences of processes $\{Z^n\}$, $\{H^n\}$, if
\[
\bigl(Y^n,H^n,Z^n \bigr) \arrowd(Y,H,Z)\qquad
\mbox{in }  \mathbb{D} \bigl(\Rp,\R^{3d}\bigr)
\]
then
\[
\bigl(\Pi \bigl(X^n \bigr),H^n,Z^n \bigr)
\arrowd(X,H,Z)\qquad  \mbox{in }  \mathbb{D} \bigl(\Rp,\R^{3d}\bigr)
\]
and
\[
\bigl(X^n_{t_1}, \ldots,X^n_{t_m},H^n,Z^n
\bigr) \arrowd (X_{t_1}, \ldots,X_{t_m},H,Z)\quad  \mbox{in }
\R^m \times\mathbb{D} \bigl(\Rp,\R^{2d}\bigr)
\]
for any $m\in\N$, any $t_1, \ldots,t_m \in\Rp$ such that $\p
(|\Delta Y_{t_i}|=0)=1, i=1, \ldots,m$, where $X=Y+K$ is a
solution of the Skorokhod problem associated with a process $Y$.
\end{longlist}
\end{corollary}

\begin{pf}
In the proof, it suffices to make the observation
that (\ref{eq3.12}) is equivalent to the fact that $
\bar{\omega}''_{( X^{n}, Z^{n})}(\delta_{n} , ) \mathop{\longrightarrow}\limits_{\mathcal{P}}0$ for
every sequence $\{\delta_{n}\}$ such that $\delta_n \downarrow0$
and to combine the deterministic results given in Theorem
\ref{tw2.15}, Corollary \ref{tw3.11} and Corollary
\ref{wniosek3.14} with the Skorokhod representation theorem.
\end{pf}

\section{Penalization methods for reflecting SDEs}\label{sec4}

Let $\spa$ be a filtered probability space, let $H$ be an
$\filt$-adapted process, and let $Z$ be an $({\mathcal{F}}_t)$ adapted
semimartingale such that $H_0\in\bar
D, Z_0=0$. Now recall
that a pair
$(X, K)$ of
$({\mathcal{F}}_t)$-adapted processes is called a solution of the
reflecting SDE (\ref{eq1.1}) if $(X, K)$ is a solution of the
Skorokhod problem associated with $Y$ defined by
\[
Y_t=H_t+\int_0^tf(X_{s-})
 \,\mathrm{d}Z_s,\qquad  t\in\Rp.
\]
We say that the SDE (\ref{eq1.1}) has a weak solution if there
exists a probability space $\spawh$ and $\filtwh$-adapted
processes $\wh H$,$\wh Z$ and $(\wh X,\wh K)$ such that $\mathcal{
L}(\wh H,\wh Z)=\mathcal{L}(H,Z)$ and $(\wh X,\wh K)$ is a solution
of the Skorokhod problem associated with $\wh Y_t=\wh
H_t+\int_0^tf(\wh X_{s-})\, \mathrm{d}\wh Z_s$, $t\in\Rp$. If any two weak
solutions $(\wh X,\wh K)$, $(\wh X',\wh K')$ of the SDE
(\ref{eq1.1}), possibly defined on two different probability
spaces, are such that $\mathcal{L}(\wh X,\wh K)=\mathcal{L}(\wh X',\wh
K')$, we say that the weak uniqueness for the SDE (\ref{eq1.1})
holds.

In this section we prove general results on weak and strong
approximations of $X$. We begin with two technical lemma. In the
first lemma, which is a simple consequence of Corollary \ref{ap4}
(\hyperref[app]{Appendix}), the sequence $\{X^n\}$ need not consist of solutions
of penalized equations, and the process $X$ need not be a solution
of (1.1).

\begin{lemma}
\label{wnJakub1}
Let $\{ X^n \}$ be a sequence of $\fn$ adapted
processes, and let $\{Z^n\}$ be
a sequence of $\fn$-adapted semimartingales satisfying \textup{(UT)} such
that
$Z^n_0=0$, $n \in\N$, $\{H^n\}$ be a sequence of processes. If
$\{ X^n \}$ is $S$-tight and there exist processes $X, H, and Z$ such
that
\[
\bigl(X^n_{t_1}, \ldots, X^n_{t_m}
,H^n,Z^n \bigr) \arrowd(X_{t_1}, \ldots,
X_{t_m},H,Z)\qquad  \mbox{in }   \R^{md} \times{\mathbb{D}} \bigl(\Rp,\R^{2d}\bigr)
\]
for any $m\in\N$ and any $  t_1, t_2,
\ldots, t_m $ from a dense subset $\mathbb Q$ of $\Rp$ and
\begin{equation}
\label{omegawar7}
\lim_{\delta\downarrow0} \limsup_{n
\rightarrow\infty} P \bigl(\bar{
\omega}''_{(X^n,Z^n)}(\delta, q)> \varepsilon
\bigr)=0, \qquad \varepsilon>0,   q \in\Rp,
\end{equation}
then for every continuous function $g \dvtx \Rd\lra\Rd\otimes\Rd$
\[
\biggl(X^n_{t_1}, \ldots, X^n_{t_m},H^n,
Z^n ,\int_0^{\cdot} \bigl\langle g
\bigl( X^n_{s-} \bigr) , \mathrm{d}Z^n_s
\bigr\rangle \biggr) \arrowd \biggl(X_{t_1},\ldots, X_{t_m},
H, Z , \int_0^{\cdot} \bigl\langle
g(X_{s-}) , \mathrm{d}Z_s \bigr\rangle \biggr)
\]
in
$\R^{md} \times\mathbb{D} (\Rp,\R^{3})$ for any $m\in\N$ and
$  t_1, t_2, \ldots,
t_m \in\mathbb Q$.
\end{lemma}

\begin{pf}
Clearly,
\begin{eqnarray}
\label{wstepRDD}
&&\bigl(g \bigl(X^n_{t_1} \bigr), \ldots, g
\bigl(X^n_{t_m} \bigr) ,X^n_{t_1},
\ldots, X^n_{t_m} ,H^n,Z^n \bigr)
\nonumber\\[-8pt]\\[-8pt]
&&\quad\arrowd \bigl(g(X_{t_1}), \ldots, g(X_{t_m}),X_{t_1},
\ldots, X_{t_m},H,Z \bigr)\qquad \mbox{in } \R^{2md} \times{\mathbb{D}} \bigl(\Rp,\R^{2d}\bigr)\nonumber
\end{eqnarray}
for any $m\in\N$ and any $  t_1, t_2,
\ldots, t_m \in\mathbb Q$. Because $\{ X^n \}$ is $S$-tight, it
follows that
$\{ g(X^n) \}$ is $S$-tight.
Similarly,
(\ref{omegawar7}) implies that
\begin{equation}
\label{wstepflomega}
\lim_{\delta\downarrow0} \limsup_{n \rightarrow\infty
}\p \bigl( \bar{
\omega}''_{(g(X^n),Z^n)}(\delta, T) > \varepsilon
\bigr)=0,\qquad  \varepsilon > 0.
\end{equation}
Combining (\ref{wstepRDD}), (\ref{wstepflomega}) and putting
$Y^n=g(X^n)$, $l=m$, $K^n_i=X^n_{t_i}$, $i=1, \ldots, m$, in
Corollary~\ref{ap4} (\hyperref[app]{Appendix}) we complete the proof.
\end{pf}

\begin{lemma}
\label{thm3.4} Let $\{ H^n \}$ be a sequence of $\fn$ adapted
processes, $H^n_0\in\bar{D}$, $n \in\N$, and let $\{Z^n\}$ be
a sequence of $\fn$-adapted semimartingales satisfying \textup{(UT)},
$Z^n_0=0$, $n \in\N$. Let $\{X^n\}$ be a sequence of solutions of
the SDE (\ref{eq1.3}). If $f\dvtx\R^d \rightarrow\Rd\otimes\Rd$
satisfies (\ref{eq1.2}) and $\{ H^n \}$ is tight in $\D$, then
$\{Y^n=X^n-H^n\}$ satisfies \textup{(UT)}.
\end{lemma}

\begin{pf}
First, we show that for every $q\in\Rp$,
\begin{equation}
\label{eq3.6a} \Bigl\{\sup_{t\leq q}\big|X^n_t\big| \Bigr\}
 \mbox{ is bounded in probability}.
\end{equation}
Let $\wh X^n$ denote the solution of the equation with
penalization term (\ref{eq3.2}) with $\wh{Y}^n=H^n$, $n\in\N$.
Because $\{ H^n \}$ is tight in $\D$, it follows by Corollary
\ref{thm2.7} that for every $q\in\Rp$, $ \{\sup_{t\leq
q}|\wh{X}^n_t|\}$ is bounded in probability.
On the other hand, $\{Z^n\}$ satisfies (UT), and thus we may and will
assume that
$Z^n_t=M^n_t+V^n_t$ and $M^n_0=V^n_0=0$, where $\{[M^n]_q\}$,
$\{|V^n|_q\}$ are bounded in probability and $|\Delta M^n|\leq c$
for some $c>0$. In this case, $\{\langle  M^n\rangle _q\}$ is bounded in
probability as well. Define $ \tau_k^n=\inf\{t; |\wh X^n_t|\vee
|V^{n}|_t\vee [M^{n}]_t\vee  \langle M^{n}\rangle_t>k\}\wedge
k, $
$n,k\in\N$. It is clear that
\begin{equation}
\label{eq4.15} \lim_{k\rightarrow+\infty}\limsup_{n\rightarrow+\infty}P \bigl(
\tau^n_k\leq q \bigr)=0, \qquad q\in\Rp.
\end{equation}
By (\ref{eq3.13}) with $p=1$ and by (\ref{eq1.2})
for every stopping time $\sigma^n$,
\begin{eqnarray*}
&&E\sup_{t<\sigma^n\wedge\tau^n_k}\big|X^{n}_t-\wh X^n_t\big|^{2}
\\
&&\quad \leq C(1) \biggl\{ E\int_0^{(\sigma^n\wedge\tau^n_k)-}\big\Vert f
\bigl(X^n_{s-} \bigr)\big\Vert ^2 \,\mathrm{d}
\bigl[M^n \bigr]_s +E\int_0^{(\sigma^n\wedge\tau^n_k)-}\big\Vert f
\bigl(X^n_{s-} \bigr)\big\Vert ^2 \,\mathrm{d}\big\langle M^n\big\rangle _s
\\
& &\hspace*{25pt}\qquad {}   +kE\int_0^{(\sigma^n\wedge\tau
^n_k)-}\big\Vert f
\bigl(X^n_{s-} \bigr)\big\Vert ^2 \,\mathrm{d}\big|V^n\big|_s
\biggr\}
\\
&&\quad \leq C(k,L) \biggl\{1+ E\int_0^{(\sigma^n\wedge\tau^n_k)-}
\sup_{u\leq s}\big|X^{n}_{u-}-\wh X^n_{u-}\big|^{2}
 \,\mathrm{d} \bigl(\big|V^n\big|+ \bigl[M^n \bigr]+\big\langle M^n\big\rangle
\bigr)_s \biggr\}.
\end{eqnarray*}
%
Therefore, for every stopping time $\sigma^n$,
\begin{eqnarray*}
& &E \sup_{t<\sigma^n} \big|X^{n,\tau^n_k-}_t-\wh
X^{n,\tau^n_k-}_t\big|^{2}
\\
&& \quad \leq C(k,L) \biggl\{ 1+ E\int_0^{\sigma^n-}
\sup_{u\leq
s}\big|X^{n,\tau^n_k-}_{u-}-\wh X^{n,\tau^n_k-}_{u-}\big|^{2}
\,\mathrm{d} \bigl(\big|V^{n,\tau^n_k-}\big|+ \bigl[M^{n,\tau^n_k-} \bigr]+\big\langle M^{n,\tau^n_k-}\big\rangle
\bigr)_s \biggr\}.
\end{eqnarray*}
Thus, by Gronwall's lemma (see, e.g., \cite{s1}, Lemma 3),
\[
E\sup_{t<\tau^n_k}\big|X^{n}_t-\wh X^n_t\big|^{2}
\leq C(k,L)\exp \bigl\{ 3k C(k,L) \bigr\},
\]
and using (\ref{eq4.15}) gives (\ref{eq3.6a}).
Combining (\ref{eq3.6a}) with (\ref{eq1.2}) shows that $\{ \sup_{t
\leq q
}\Vert f(X^n_{t-})\Vert  \}$ is also bounded in probability, $q\in\Rp$.
Thus, the sequence of stochastic integrals
$\{\int_0^{\cdot}\langle f(X^n_{s-}) , \mathrm{d}Z^n_s \rangle\}$
satisfies (UT). Because of Corollary \ref{thm2.7}, the proof is
complete.
\end{pf}

We can now formulate our main theorem.

\begin{theorem}
\label{tw3}
Let $\{ H^n \}$ be a sequence of $\fn$ adapted
processes, $H^n_0\in\bar{D}$, $n \in\N$, and let $\{Z^n\}$ be
a sequence of $\fn$ adapted semimartingales satisfying \textup{(UT)},
$Z^n_0=0$, $n \in\N$. Let $\{X^n\}$ be a sequence of solutions of
the SDE (\ref{eq1.3}). If $f\dvtx\R^d \rightarrow\Rd\otimes\Rd$ is
continuous, satisfies (\ref{eq1.2}), and ${
(H^n,Z^n)\arrowd(H,Z)}$ in $\DD$, then, for every continuous
function $g \dvtx \Rd\lra\Rd\otimes\Rd$,
\begin{longlist}
\item[(i)] $\{(X^n_{t_1}  ,  \ldots  ,  X^n_{t_m}  ,  H^n
,  Z^n  ,
\int_0^{\cdot}\langle g(X^n_{s-})  ,   \mathrm{d}Z^n_s \rangle)\}$ is
tight in $\R^{md}  \times  \Diii $ and its every limit point
has the form
\[
\biggl(X_{t_1},\ldots,X_{t_m},H,Z,\int_0^{\cdot}
\bigl\langle g(X_{s-}),\mathrm{d}Z_s \bigr\rangle \biggr)
\]
for any $m\in\N$ and any
$t_1,\ldots ,t_m \in\Rp$ such that $\p(|\Delta Z_{t_i}|=0)= \p
(|\Delta H_{t_i}|=0)=1$, $i=1,\ldots ,m$, where $X$ is weak solution
of the SDE (\ref{eq1.1}),
\item[(ii)]
if (\ref{eq1.1}) has a unique weak solution $X$, then
\[
\biggl(X^n_{t_1},\ldots ,X^n_{t_m},H^n,Z^n,
\int_0^{\cdot}  \bigl\langle g
\bigl(X^n_{s-} \bigr),\mathrm{d}Z^n_s \bigr
\rangle \biggr) \arrowd \biggl(X_{t_1},\ldots ,X_{t_m},H,Z,\int
_0^{\cdot}  \bigl\langle g(X_{s-}),\mathrm{d}Z_s
\bigr\rangle \biggr)
\]
in $\R^{md}\times\Diii$, for any $m\in\N$ and any $t_1,\ldots ,t_m
\in\Rp$ such that $\p(|\Delta Z_{t_i}|=0)= \p(|\Delta
H_{t_i}|=0)=1$, $i=1,\ldots ,m$,
\item[(iii)] if $X$ has a unique weak
solution, then $ {X^n\arrowds X}$.
\end{longlist}
\end{theorem}

\begin{pf}
(i) First, note that by Lemma \ref{thm3.4}, $X^n$ has the
form
\begin{equation}
\label{pkt1} X^n_t=H^n_t+W^n_t,
 \qquad t \in\Rp,
\end{equation}
where $\{ W^n \}$ is a sequence of semimartingales satisfying
(UT). We set $\gamma^{i}_0=0$, $\gamma^{i}_{k+1}=
\min(\gamma^{i}_{k}+\delta^i_k, \inf\{ t > \gamma^{i}_{k};
|\Delta H_t | > \delta^i \})$ and $\gamma^{ni}_0=0$, $\gamma^{ni}_{k+1}= \min(\gamma^{ni}_{k}+\delta^i_k, \inf\{ t > \gamma^{ni}_{k}; |\Delta H^n_t | > \delta^i \})$, where $\{ \delta^i
\}$, $\{ \{ \delta^i_k \} \}$ are families of positive constants
such that $\delta^i \downarrow0$, $\delta^i/2 \leq\delta^i_k
\leq\delta^i$ $P(|\Delta H_t|=\delta^i, t \in\Rp)=0$,
$P(|\Delta H_{\gamma_k^{i}+\delta^i_k}| =0 )=1$. For every
$i\in\N$, define a new sequence, $\{H^{ni}\}$, of processes by putting
$H^{ni}_t=H^n_{\gamma^{ni}_{k}}$, $t \in[\gamma^{ni}_{k},
\gamma^{ni}_{k+1} )$, $k \in\No$, $n\in\N.$ Then, using the
continuous mapping theorem, we have $(H^n,H^{ni}) \arrowd
(H,H^i)$ in $\DD$, which implies that
\begin{equation}
\label{ni} \lim_{i \rightarrow
\infty}\limsup_{n \rightarrow\infty} \p \Bigl(
\sup_{t \leq q}\big|H^{ni}_t-H^n_t\big|
\geq\varepsilon \Bigr)=0, \qquad  \varepsilon>0,  q \in\Rp.
\end{equation}
If $ X^{ni}$ denotes the solution of the equation
\begin{eqnarray*}
X^{ni}_t & = & H^{ni}_t+ \int
_0^{t} \bigl\langle f \bigl(X^n_{s-}
\bigr) , \mathrm{d}Z^n_s \bigr\rangle+ K^{ni}_t,
\qquad  t \in\Rp, n,i \in\N,
\end{eqnarray*}
then, by (\ref{ni}) and Corollary \ref{cor2.8}, $ \lim_{i
\rightarrow\infty}\limsup_{n \rightarrow\infty}
\p(\sup_{t \leq q}|X^{ni}_t-X^n_t| \geq\varepsilon)=0$,
$\varepsilon>0$, \mbox{$q\in\Rp$}.
Furthermore, it is well
known that for continuous $f\dvtx \Rd\lra\Rd\otimes\Rd$, one can
construct a sequence $\{ f^i \}$ of functions such that $f^i \in
\cdw$, $i\in\N$ and $\sup_{x \in{K}}\Vert f^i(x) - f(x) \Vert  \lra0$
for any compact subset $K \subset\Rd$. If we set
\begin{eqnarray*}
Y^{ni}_t & = & H^{n}_t+ \int
_0^{t} \bigl\langle f^i
\bigl(X^{ni}_{s-} \bigr) , \mathrm{d}Z^n_s
\bigr\rangle, \qquad t \in\Rp, n \in\N,
\\
Y^{n}_t & = & H^{n}_t+ \int
_0^{t} \bigl\langle f \bigl(X^{n}_{s-}
\bigr) , \mathrm{d}Z^n_s \bigr\rangle, \qquad  t \in\Rp, n \in\N,
\end{eqnarray*}
then $ \lim_{i \rightarrow\infty}\limsup_{n \rightarrow\infty
}
\p(\sup_{t \leq q}|Y^{ni}_t-Y^n_t| \geq\varepsilon)=0$,
$\varepsilon
> 0$, $q\in\Rp.$
Because $f^i \in\cdw$ and $\{X^{ni}\}$ satisfies (UT), $\{
f^i(X^{ni}) \}$ satisfies (UT) as well. By \cite{s2}, Lemma 4.3,
the sequences $\{ Y^{ni} \}$ and $\{ Y^n \}$ are tight in $\D$.
Moreover, we can see that $\{ (Y^n,H^n,Z^n) \}$ is tight in
$\DDD$. Because $X^n$ is a solution of the equation with
a penalization term of the form (\ref{eq1.5}), it follows by
Corollary \ref{thm2.7}, Proposition \ref{ap15} (\hyperref[app]{Appendix}), and
Corollary \ref{wnzb}(i) that
\begin{equation}
\label{RZZ} \bigl\{ X^n \bigr\} \mbox{ is $S$-tight}
\end{equation}
and
\begin{equation}
\label{omega} \lim_{\delta\downarrow0} \limsup_{n \rightarrow\infty
}\p \bigl( \bar{
\omega}''_{(X^n,Z^n)}(\delta, q) > \varepsilon
\bigr)=0, \qquad \varepsilon > 0, q\in\Rp.
\end{equation}
Assume that there exists a subsequence $\{ n' \} \subset\{ n \}$
such that $(Y^{n'}, H^{n'} , Z^{n'} ) \arrowd(\wh Y, \wh H , \wh
Z )$ in $\DDD$, where $\mathcal{L}(\wh H,\wh Z)=\mathcal{L}(H,Z)$.
Then, by Corollary \ref{wnzb}(ii),
\begin{equation}
\label{1RDD} \bigl(X^{n'}_{t_1}, \ldots ,X^{n'}_{t_m}
, H^{n'}, Z^{n'} \bigr) \arrowd(\wh X_{t_1},
\ldots, \wh X_{t_m}, \wh H, \wh Z) \qquad \mbox{in } \R^{dm}
\times{\mathbb{D}} \bigl(\Rp,\R^{2d}\bigr)
\end{equation}
for any $m\in\N$ and any $t_1, \ldots,t_m \in\Rp$ such that
$\wh\p( |\Delta\wh Y_{t_i}|=0)=1$, $i=1, \ldots,m $, where
$\wh X$ is a solution of the Skorokhod problem associated with
$\wh Y$. Combining (\ref{RZZ})--(\ref{1RDD}) with Lemma
\ref{wnJakub1} yields
\[
\biggl(X^{n'}_{t_1}, \ldots, X^{n'}_{t_m},
H^{n'} , Z^{n'} , \int_0^{\cdot}
\bigl\langle g \bigl(X^{n'}_{s-} \bigr), \mathrm{d}Z^{n'}_s
\bigr\rangle \biggr) \arrowd \biggl(\wh X_{t_1}, \ldots, \wh
X_{t_m},\wh H , \wh Z , \int_0^{\cdot}
\bigl\langle g(\wh X_{s-}) , \mathrm{d}\wh Z_s \bigr\rangle \biggr)
\]
in $\R^{md} \times\DDD$ for any $m\in\N$
and any $t_1, \ldots,t_m \in\Rp$ such that $\wh\p( |\Delta\wh
Y_{t_i}|=0)=1$, $i=1, \ldots,m $. Thus, in particular, putting
$g=f$, we obtain that $\wh Y = \wh H + \int_0^{\cdot} f(\wh
X_{s-}) \,\mathrm{d}\wh Z_s$, which implies that $\wh X$ is a weak solution
of the SDE (\ref{eq1.1}). Because $\{t;\wh\p( |\Delta\wh
Y_{t}|=0)=1\}\subset\{t;\wh\p( |\Delta\wh H_{t}|=0)=1  \mbox{ and } \wh\p( |\Delta\wh Z_{t}|=0)=1\}$, the proof of (i) is
complete.
(ii) Follows immediately from (i).

(iii) Because $\{X^n\}$ is $S$-tight, the desired result follows
from the convergence of finite-dimensional distributions of $X^n$ to
these of $X$ proven in part (ii) and from Corollary \ref{ap2}
(\hyperref[app]{Appendix}).
\end{pf}



We now consider an array $\{\{\tk\}\}$ of nonnegative numbers
such that the $n$th row $T_n=\{\tk\}$ forms a
partition of $\Rp$ such that $0=t_{n,0}<t_{n,1}<\cdots ,$
$\lim_{k\rightarrow\infty}\tk=+\infty$ and $
\max_k (\tk-\tkk)\longrightarrow0 \mbox{
as }  n\rightarrow+\infty$.  For the array $\{\{\tk\}\}$, we
define a sequence of summation rules $\{\rnn\}, \rnn\dvtx\Rp\lra\Rp$
by $\rnn_t=\max\{\tk;\tk\leq t\}$, and then for a fixed adapted
process $H$ and a semimartingale $Z$, we define $
H\rn=H_{\rho^n_\cdot}$, $ Z\rn=Z_{\rho^n_\cdot}$, that is,
\[
H\rn_t=H_{\tk},\qquad  Z\rn_t=Z_{\tk}
 \qquad\mbox{for }   t\in[\tk,t_{n,k+1}), k\in\N\cup\{0\}, n\in\N.
\]

Let $\{\ox\}$ be a sequence of solutions to equations with
penalization terms
driven by $\{Z\rn\}$, that is,
\begin{equation}
\label{eq3.24} \ox_t=H\rn_t+\int_0^tf\bigl(
\ox_{s-}\bigr) \,\mathrm{d}Z_s\rn-n\int_0^t
\bigl(\ox_s-\Pi \bigl(\ox_s\bigr) \bigr) \,\mathrm{d}s, \qquad  t\in\Rp, n
\in \N.
\end{equation}
The special form of $Z\rn$ implies that
\[
   \bar X^n_t=\cases{ %
 H_0, &\quad$t=0$,
\vspace*{2pt}\cr
\Pi \bigl(\bar X^n_{t_{n,k}} \bigr)+ \bigl( \bar
X^n_{t_{n,k}}-\Pi \bigl( \bar X^n_{t_{n,k}}
\bigr) \bigr)\mathrm{e}^{-n(t-t_{n,k})}, & \quad$t \in (t_{n,k},t_{n,k+1}), k \in
\N\cup\{0\}$,
\vspace*{2pt}\cr
\bar X^n_{(t_{n,k+1})-}+( H_{t_{n,k+1}}- H_{t_{n,k}})&
\cr
 \quad+f \bigl(\bar X^n_{(t_{n,k+1})-} \bigr) ( Z_{t_{n,k+1}}-
Z_{t_{n,k}}), &\quad $t=t_{n,k+1}, k \in\N\cup\{0\}$. }
\]

\begin{corollary}
\label{cor4.2} Let $\{\ox\}$ be a sequence of solutions of
(\ref{eq3.24}). If $f$ is continuous and satisfies (\ref{eq1.2}),
then
\begin{longlist}
\item[(i)] $\{ \bar X^n\}$ is $S$-tight, and its every limit
point $X$ is a weak solution of the SDE (\ref{eq1.1}),
\item[(ii)] if the SDE (\ref{eq1.1}) has a unique weak solution,
then $ {\bar X^n\arrowds X}$.
\end{longlist}
\end{corollary}

\begin{pf}
It is easily seen that $(H\rn,Z\rn)\rightarrow(H,Z)$ almost
surely in $\Djj$. On the other hand, based on the theorem of
Bichteler, Dellacherie and Mokobodzki, the sequence $\{Z\rn\}$
of discrete semimartingales
satisfies (UT). Therefore, the assumptions of Theorem \ref{tw3}
are satisfied. Thus, (i) follows from (\ref{RZZ}), whereas (ii)
follows from
Theorem \ref{tw3}(iii).
\end{pf}

We note that Corollary \ref{cor4.2}(i) implies the existence
of a weak solution of the SDE (\ref{eq1.1}).

In the sequel, we consider convergence in probability of solutions
of equations with penalization terms. We assume that the SDE
(\ref{eq1.1}) has the pathwise uniqueness property; that is, for any\vspace*{1pt}
two solutions $\wh X$, $\wh X'$ of the SDE (\ref{eq1.1})
corresponding to processes $(\wh H,\wh Z)$, $(\wh H',\wh Z')$ and
defined on a probability space $\spawh$ with filtration $\filtwh$,
the following implication holds:
\[
\wh P \bigl((\wh H_t,\wh Z_t)= \bigl(\wh
H'_t,\wh Z'_t \bigr);t\in\Rp
\bigr)=1\quad\Rightarrow\quad \wh P \bigl(\wh X_t=\wh X'_t;t
\in\Rp \bigr)=1.
\]
It is well known that
the existence of weak solutions and the pathwise uniqueness property
implies the existence of a unique strong solution on arbitrary
probability space $\spa$ with given adapted process $H$ and the
semimartingale $Z$ (see, e.g., Yamada and Watanabe \cite{jw} and
Jacod and M\'emin \cite{jm} for the case of general
semimartingales). The classical example of equation with the
pathwise uniqueness property and non-Lipschitz coefficient $f$
was given by Yamada and Watanabe \cite{jw} in the case of
diffusion equations. Tudor \cite{tu} proved that
this example works for SDEs driven by general
semimartingales. Using, for example, \cite{s2}, Lemma C.3, we can
give the following version of his result:

\begin{example}
Assume that $f$ is continuous, satisfies (\ref{eq1.2}) and
\[
\big\Vert f(x)-f(y)\big\Vert ^2\leq\rho \bigl(|x-y|^2 \bigr),\qquad   x,y
\in \bar D,
\]
where $\rho\dvtx\Rp\to\Rp$ is strictly increasing and concave,
$\rho(0)=0$ and $\int_{0+}\frac{\mathrm{d}u}{\rho(u)}=+\infty$. Then the
SDE (\ref{eq1.1}) has the pathwise uniqueness property.
\end{example}

\begin{corollary}
\label{cor4.3} Moreover, under the assumptions of Theorem \ref{tw3}, if
$(H^n,Z^n)\arrowp(H,Z)$ in $\DD$ and the SDE
(\ref{eq1.1}) has the pathwise uniqueness property, then for any
continuous function $g \dvtx \Rd\lra\Rd\otimes\Rd$,
\begin{longlist}
\item[(i)]
\[
 {
\biggl(X^n_{t_1},\ldots ,X^n_{t_m},H^n,Z^n,\int_0^{\cdot} \big\langle
g\bigl(X^n_{s-}\bigr),\mathrm{d}Z^n_s\big\rangle \biggr) \arrowp
\biggl(X_{t_1},\ldots ,X_{t_m},H,Z,\int_0^{\cdot} \big\langle
g(X_{s-}),\mathrm{d}Z_s\big\rangle \biggr)}
\]
in $\R^{md}\times\Diii$, for
any
$m\in\N$ and any $t_1,\ldots ,t_m \in\Rp$ such that $\p(|\Delta
Z_{t_i}|=0)= \p(|\Delta H_{t_i}|=0)=1$, $i=1,\ldots ,m$,
\item[(ii)] $ {X^n\arrowdps
X}$,
\end{longlist}
where $X$ is a unique strong solution of the SDE (\ref{eq1.1}).
\end{corollary}

\begin{pf}
In the proof, it suffices to use Theorem \ref{tw3} and
repeat arguments from
the proof of Theorem 1(ii) in \cite{s1}. Fix $B\in{\mathcal{F}}$,
$\p(B)>0$ and define ${\mathcal{Q}}_B=\p(A|B)$ for every $A\in\mathcal{
F}$. Obviously, ${\mathcal{Q}}_B\ll\p$ and $(\frac{\mathrm{d}\mathcal{
Q}}{\mathrm{d}\p})=\frac{\mathbf{1}_B}{\p(B)}$. Let $\{H^n\}$, $\{Z^n\}$
be sequences of processes satisfying the assumptions of Corollary
\ref{cor4.3}. Then $\{Z^n\}$ is a sequence of semimartingales on
$(\Omega,{\mathcal{F}},{\mathcal{Q}}_B)$ satisfying (UT) (see, e.g.,
\cite{s1}, Lemma 4) and $(H^n,Z^n)\arrowdq(H,Z)$ in $\Dii$.
Moreover, the stochastic integral $\int f(X^{n}_{s-})\, \mathrm{d}Z^n_s$,
calculated with respect to $\p$ is, for ${\mathcal{Q}}_B$, almost all
$\omega\in\Omega$ equal to the integral calculated with respect to
${\mathcal{Q}}_B$ and
\[
X^{n}_t=H^n_t+\int
_0^tf \bigl(X^n_{s-}
\bigr)  \,\mathrm{d}Z^n_s-n\int_0^t\bigl(X^n_s-
\Pi \bigl(X^n_s \bigr)\bigr)\, \mathrm{d}s, \qquad  t\in\Rp,  \mathcal{
Q}_{B}\mbox{-a.e.}, n\in\N.
\]
Consequently, by Theorem \ref{tw3}(ii), $
( X^{n}_{_1},\ldots ,X^n_{t_m}, H^n , Z^n, Y^n )\arrowdq( X , H
, Z, Y )$,
in $\R^m\times\Diii$,
where $Y^n=\int_0^{\cdot}g(X^n_{s-}) \,\mathrm{d}Z^n_s$, $Y=\int_0^{\cdot}g
(X_{s-}) \,\mathrm{d}Z_s$. Thus, for all bounded and continuous mappings
$\Phi , \Phi\dvtx\R^m\times\Dii\lra\R$,
\[
\lim_{n\rightarrow+\infty}\int_{\Omega}\Phi \bigl(X^{n}_{t_1},\ldots ,X^n_{t_m},H^n,Z^n,Y^n
\bigr) \,\mathrm{d}{\mathcal{Q}}_{B} =\int_{\Omega}
\Phi(X_{t_1},\ldots ,X_{t_m},H,Z,Y) \,\mathrm{d}{\mathcal{Q}}_{B} ,
\]
or, equivalently,
\begin{equation}
\label{eq3.22} \lim_{n\rightarrow+\infty}\int_{B}\Phi
\bigl(X^{n}_{t_1},\ldots ,X^n_{t_m}
,H^n,Z^n,Y^n \bigr) \,\mathrm{d}{\p} =\int
_{B}\Phi(X_{t_1},\ldots, X_{t_m},H,Z,Y) \,\mathrm{d}{
\p} . 
\end{equation}
Because (\ref{eq3.22}) holds for all $B\in{\mathcal{F}} , \mathcal{
P}(B)>0$ and all bounded continuous mappings
$\Phi\dvtx\R^m\times\Diii\lra\R$, the proof of (i) is complete. Using
(i), the claim (ii) follows readily.
\end{pf}

\begin{corollary}
\label{cor4.4} Assume that $f$ is continuous and satisfies
(\ref{eq1.2}), and that (\ref{eq1.1}) has the pathwise uniqueness
property. Let $X$ be a strong solution of (\ref{eq1.1}), and let
$\{\ox\}$ be a sequence of solutions of (\ref{eq3.24}). Then
\begin{longlist}
\item[(i)] ${ \bar X^n_{t}\arrowp X_t}$
for every $t\in\Rp$ such that $P(\Delta H_{t}=0)=P(\Delta
Z_{t}=0)=1$,
\item[(ii)]
for any continuous $g\dvtx\R^d
\rightarrow\Rd\otimes\Rd$,
\[
\sup_{t\leq q, t\in T_n}\bigg|\int_0^{t}g \bigl(\bar
X^n_{s-} \bigr) \,\mathrm{d}Z\rn_s -\int
_0^{t}g(X_{s-}) \,\mathrm{d}Z_s\bigg|
\arrowp0, \qquad q\in\Rp,
\]
\item[(iii)] ${ \bar X^n\arrowdps
X \mbox{ in } \D,}$
\item[(iv)] if, moreover, $H,Z$
are processes with continuous trajectories, then
\[
\sup_{t\leq q}\big|\bar X^n_t-X_t\big|\arrowp0,
 \qquad q\in\Rp.
\]
\end{longlist}
\end{corollary}

\begin{pf}
We can see that $(H\rn,Z\rn)\rightarrow(H,Z)$ almost surely
in $\Dii$, and, given the theorem of Bichteler, Dellacherie, and
Mokobodzki, the sequence $\{Z\rn\}$ of discrete semimartingales
satisfies (UT). Therefore, by
Corollary \ref{cor4.3}, the conclusions (i) and (iii) follow.
Moreover, by Corollary \ref{cor4.3}(i), for every continuous $g,$
\[
\biggl(\int_0^{\cdot}g \bigl(\bar
X^n_{s-} \bigr) \,\mathrm{d}Z\rn_s,Z\rn \biggr) \arrowp
\biggl(\int_0^{\cdot}g(X_{s-})
 \,\mathrm{d}Z_s,Z \biggr),\qquad \mbox{in }  \R^m\times \mathbb{D} \bigl(\Rp,\R^{2d}\bigr).
\]
Using this and arguing as in the proof of Theorem 3 in
\cite{s1}, the conclusion (ii) also follows.

To prove (iv), we first observe that by Corollary
\ref{cor4.3}(i),
\[
Y^n=H\rn+\int_0^{\cdot}f \bigl(\bar
X^n_{s-} \bigr)\, \mathrm{d}Z\rn_s\arrowp H+\int
_0^{\cdot}f(X_{s-}) \,\mathrm{d}Z_s=Y
 \qquad\mbox{in }  {\mathbb{D}} \bigl(\Rp,\Rd\bigr),
\]
where
$Y$ is a process with continuous trajectories. From this, and by
Corollary \ref{tw3.11}(iii), the result follows readily.
\end{pf}

\begin{appendix}\label{app}

\renewcommand{\theequation}{\textup{\Alph{section}}.\arabic{equation}}
\setcounter{equation}{0}

\renewcommand{\theremark}{\Alph{section}.\arabic{remark}}
\setcounter{remark}{0}

\section*{Appendix: The topology $S$}

The $S$ topology on the space $\D$ of $\Rd$-valued functions
that are right-continuous and have left-hand limits was
introduced by Jakubowski \cite{ja}. It
is weaker than the Skorokhod topology $J_1$
but stronger than the Meyer--Zheng topology considered in
\cite{mz}. We collect here only basic
definitions and properties of the $S$ topology; more details can be
found in Jakubowski \mbox{\cite{ja0,ja}}.

\begin{proposition}
\textup{(i)} $K\subset\D$ is relatively $S$-compact
if and only if
\begin{equation}
\label{eq5.1} \sup_{x\in K}\sup_{t\leq q}|x_t|<+\infty,
\qquad  q\in\Rp
\end{equation}
and for all $a<b$, $a,b\in\R$
\begin{equation}
\label{eq5.2} \sup_{x\in K}N^{a,b}(x,q)<+\infty,\qquad  q\in\Rp,
\end{equation}
where $N^{a,b}$ is the usual number of up-crossings given levels
$a<b$, that is, $N^{a,b}(x,q)\geq k$ if one can find numbers $0\leq
t_1<t_2<\cdots <t_{2k-1}<t_{2k}\leq q$ such that $x_{t_{2i-1}}<a$ and
$x_{t_{2i}}>b$, $i=1,2,\ldots ,k$.

\textup{(ii)} $x^n$ converges
to $x$ in the $S$-topology if and only if $\{x^n\}$ satisfies
(\ref{eq5.1}), (\ref{eq5.2}) and
in
every subsequence $\{n_k\}$, one can find a further subsequence
$\{n_{k_l}\}$ and a
dense subset $\mathbb Q\subset
\Rp$ such that $x^{n_{k_l}}_t\rightarrow x_t$, $t\in\mathbb Q$.
\end{proposition}

\begin{corollary}
\label{ap1} If $\{x^n\}$ is relatively $S$-compact and there
exists a dense subset $\mathbb Q$ such that for every $t\in
\mathbb Q$, $x^n_t\rightarrow x_t$, then $\{x^n\}$ converges to
$x$.
\end{corollary}

Recall that the
sequence of processes $\{X^n\}$ converges weakly to $X$ in the
$S$ topology ($ {X^n\arrowds X}$) if in every
subsequence $\{X^{n_{k}}\}$, we can find a further subsequence
$\{X^{n_{k_l}}\}$ and stochastic processes $\{Y_l\}$ defined on
$([0,1],\mathcal{ B}_{[0,1]},{{\ell}})$, such that the laws of $Y_l$
and $X^{n_{k_l}}$ are the same, $l\in\N$, for each
$\omega\in[0,1]$
$Y_l(\omega)$ converges to $Y(\omega)$ in the $S$ topology, and
for each $\varepsilon>0$, there exists an $S$-compact subset
$K_{\varepsilon}\subset\D$ such that
\[
P \bigl( \bigl\{\omega\in[0,1]\dvt Y_l(\omega)\in K_{\varepsilon},l=1,2,\ldots
\bigr\} \bigr)>1-\varepsilon.
\]

\begin{proposition}
\label{ap0} The following two conditions are equivalent:
\begin{longlist}
\item[(i)] $\{X^n\}$ is $S$-tight.
\item[(ii)] $\{X^n\}$ is relatively compact with respect to
the convergence ``$ {\arrowds}$''.
\end{longlist}
\end{proposition}

\begin{proposition} \label{ap15}
Let $\{X^n\}$ be a sequence of
processes of the form $X^n=H^n+Z^n$, $n\in\N$, where $\{H^n\}$ is
tight in $\D$ and $\{Z^n\}$ is a sequence of semimartingales
satisfying \textup{(UT)}. Then $\{X^n\}$ is $S$-tight.
\end{proposition}

\begin{corollary}
\label{ap2} If $\{X^n\}$ is $S$-tight and there exists a dense
subset $\mathbb Q\subset\Rp$ such that for every $m\in\N$ and
every $t_1,t_2,\ldots, t_m\in\mathbb Q$
\[
\bigl(X^n_{t_1},X^n_{t_2},\ldots ,X^n_{t_m}
\bigr)\arrowd \bigl(X^n_{t_1},X^n_{t_2},\ldots ,X^n_{t_m}
\bigr)\qquad \mbox{in }  \Rd
\]
then $ {X^n\arrowds X}$.
\end{corollary}

\begin{theorem}\label{ap3}
Suppose $\{Z^n\}$ satisfies \textup{(UT)} and is tight
in $\D$ and $\{Y^n\}$ is $S$-tight. If there exists a dense
subset $\mathbb Q\subset\Rp$ such that for any $m\in\N$, any
$t_1,\ldots ,t_m$, $t_j\in\mathbb Q$, $j=1,\ldots ,m$,
\[
\bigl(Y^n_{t_1},Z^n_{t_1},\ldots ,Y^n_{t_m},Z^n_{t_m}
\bigr) \arrowd (Y_{t_1},Z_{t_1},\ldots ,Y_{t_m},Z_{t_m})
 \qquad \mbox{in }  \R^{2md},
\]
where both processes $Y$ and $Z$
have trajectories in $\D$ and there are no oscillations of
$Y^n$ preceding oscillations of $Z^n$ i.e.
\begin{equation}
\label{a}\lim_{\delta\to0}\limsup_{n\to\infty} P \bigl(\bar
\omega''_{(Y^n,Z^n)}(\delta,q)>\varepsilon
\bigr)=0,\qquad  \varepsilon>0, q\in\Rp
\end{equation}
then we have
$ {(\int_0^{\cdot}Y^n_{s-}\,\mathrm{d}Z^n_s, Z^n)\arrowd
(\int_0^{\cdot} Y_{s-}\,\mathrm{d}Z_s, Z) \mbox{ in }
\Db}$.
\end{theorem}

\begin{corollary}\label{ap4}
Suppose $\{Z^n\}$ satisfies \textup{(UT)} and is tight in $J_1$ and
$\{Y^n\}$ is $S$-tight. If there exist sequences of random vectors
$\{K^n_1 \}, \ldots,\{K^n_l \}$ and a dense subset $\mathbb
Q\subset\Rp$ such that for any $m\in\N$, any $t_1,\ldots ,t_m$,
$t_j\in\mathbb Q$, $j=1,\ldots ,m$,
\[
\bigl(Y^n_{t_1}, \ldots, Y^n_{t_m}
,K^n_1, \ldots, K^n_l
,H^n,Z^n \bigr) \arrowd(Y_{t_1}, \ldots,
Y_{t_m},K_1, \ldots,K_l, H,Z)
\]
in
$\R^{(m+l)d}\times\Djj$, where both processes $Y$ and $Z$ have
trajectories in $\D$ and (\ref{a}) holds true then
\[
\biggl(K^n_{1}, \ldots, K^n_{l},H^n,
Z^n ,\int_0^{\cdot} \bigl\langle
Y^n_{s-} , \mathrm{d}Z^n_s \bigr\rangle
\biggr) \arrowd \biggl(K_{1},\ldots, K_{l}, H, Z , \int
_0^{\cdot}\langle Y_{s-} ,
\mathrm{d}Z_s \rangle \biggr)
\]
in $\R^{ld} \times\Daaa$.
\end{corollary}
\end{appendix}

\section*{Acknowledgements}
We thank the referee for careful reading of the paper and many valuable
remarks.
We also thank prof. Stanis\l aw Kasjan from
Nicolaus Copernicus University for communicating us example~(\ref{eq2.1}). The second author was supported by the Polish
Ministry of Science and Higher Education Grant N N201372436.


%

\printhistory

\end{document}